\documentclass[journal,article,submit,moreauthors,pdftex,10pt,a4paper]{mdpi}

\usepackage{array}

\usepackage[utf8]{inputenc} 
\usepackage[T1]{fontenc}    
\usepackage{hyperref}       
\usepackage{url,mathtools}            
\usepackage{booktabs}       
\usepackage{amsfonts}       
\usepackage{nicefrac}       
\usepackage{microtype}      
\usepackage{graphicx}
\usepackage{color}
\usepackage{ulem}

\usepackage{algorithm,algorithmic}
\usepackage{amsthm,amsmath,amssymb}
\usepackage{natbib}
\usepackage[super]{nth}
\usepackage{enumerate}

\newcommand{\yc}[1]{\textcolor{black}{#1}}
\newcommand{\ycc}[1]{\textcolor{black}{#1}}
\newcommand{\yao}[1]{\textcolor{black}{#1}}

\firstpage{1} 
\articlenumber{x}
\doinum{10.3390/------}
\pubvolume{xx}
\pubyear{2017}
\copyrightyear{2017}
\externaleditor{Academic Editor: name}
\history{Received: date; Accepted: date; Published: date}

\preto{\abstractkeywords}{\nolinenumbers}

\Title{\yao{Sequential change-point detection via online convex optimization}}

\Author{Yang Cao$^{1}$, Liyan Xie$^{1}$, Yao Xie$^{1}$*, Huan Xu$^{1}$}

\AuthorNames{Yang Cao, Liyan Xie, Yao Xie, and Huan Xu}

\address{%
$^{1}$ \quad H. Milton Stewart School of Industrial and Systems Engineering, Georgia Institute of Technology; \{caoyang, lxie49\}@gatech.edu, \{yao.xie, huan.xu\}@isye.gatech.edu\\
}

\corres{Correspondence: yao.xie@isye.gatech.edu}

\abstract{Sequential change-point detection when the distribution parameters are unknown is a fundamental problem in statistics and machine learning. When the post-change parameters are unknown, \yao{we consider a set of detection procedures based on sequential likelihood ratios with non-anticipating estimators constructed using online convex optimization algorithms such as online mirror descent, which provides a more versatile approach to tackle complex situations where recursive maximum likelihood estimators cannot be found.}
When the underlying distributions belong to a exponential family and  the estimators satisfy the logarithm regret property, we show that this approach is nearly second-order asymptotically optimal. 
This means that the upper bound for the false alarm rate of the algorithm (measured by the average-run-length) meets the lower bound asymptotically up to a log-log factor when the threshold tends to infinity. Our proof is achieved by making a connection between sequential change-point and online convex optimization and leveraging the logarithmic regret bound property of online mirror descent algorithm. Numerical and real data examples validate our theory.
}


\begin{document}

\setcounter{section}{0} 

\section{Introduction}

Sequential analysis is a classic topic in statistics concerning {\it online} inference from a sequence of observations.  The goal is to make statistical inference {\it as quickly as possible}, while controlling the false-alarm rate. \yc{An important sequential analysis problem commonly studied is} sequential change-point detection \cite{siegmund1985sequential}. \yc{It} arise\yc{s} from various applications including online anomaly detection, statistical quality control, biosurveillance, financial arbitrage detection and network security monitoring (see, e.g., \citep{parametric_changepoint_2000,siegmund_new_survey,  tartakovsky2014sequential}).

We are interested in \yc{the sequential change-point detection problem with \textit{known} pre-change parameters but \textit{unknown} post-change parameters.} \yc{Specifically, }given a sequence of samples $X_1$, $X_2$, $\ldots$, \yc{we assume} that they are \yc{independent and identically distributed (i.i.d.)} with certain distribution $f_\theta$ parameterized by $\theta$, and the values of $\theta$ are different before and after \yc{some unknown time called }the \textit{change-point}. \yc{We further} assume that \yc{the parameters before the change-point are known.} This is reasonable since \yc{usually it is relatively easy to obtain the reference data for the normal state, so that the parameters in the normal state can be estimated with good accuracy.} After the change\yc{-point}, \yc{however, }the value\yc{s} of the parameter\yc{s} \yc{switch} to \yc{some} {\it unknown} value\yc{s}, \yc{which} represent anomalies or novelties that need to be discovered.

\subsection{Motivation: Dilemma of CUSUM and generalized likelihood ratio (GLR) statistics}

Consider change-point detection with unknown \yc{post-change }parameters. A commonly used change-point detection method is the so-called CUSUM procedure \citep{tartakovsky2014sequential} that can be derived from likelihood ratios. Assume that before the change, the samples $X_i$ follow a distribution $f_{\theta_0}$ and after the change the samples $X_i$ follow another distribution $f_{\theta_1}$. CUSUM procedure has a recursive structure: \yc{initialized} with $W_0 = 0$, the likelihood-ratio statistic can be computed according to $W_{t+1} = \max\{W_t + \log( f_{\theta_1}(X_{t+1})/f_{\theta_0}(X_{t+1})), 0\}$, and a change-point is detected whenever $W_t$ exceeds a pre-specified threshold. Due to the recursive structure, CUSUM is \yc{memory and computation efficient} since it does not need to store the historical data and only needs to record the value of $W_t$. 
\yc{The performance of CUSUM depends on the choice of the post-change parameter $\theta_1$; in particular, there must be a well-defined notion of ``distance'' between $\theta_0$ and $\theta_1$.} However, \yc{the choice of $\theta_1$} is somewhat subjective. \yc{Even if in practice a reasonable choice of $\theta_1$  is the ``smallest'' change-of-interest}, in the multi-dimensional setting, it is hard to define what the ``smallest'' change would mean. Moreover, when the assume\yc{d} parameter $\theta_1$ deviates significantly from the true parameter value, CUSUM may suffer a severe performance degradation \cite{granjon2013cusum}. 

An alternative approach is the Generalized Likelihood Ratio (GLR) statistic based \yc{procedure} \cite{basseville1993detection}. The GLR statistic finds the maximum likelihood estimate (MLE) of the post-change parameter and plugs it back to the likelihood ratio to form the detection statistic. To be more precise, for each hypothetical change-point location $k$, the corresponding post-change samples are $\{X_{k+1}, \ldots, X_t\}$. Using these samples, one can form the MLE denoted as $\hat{\theta}_{k+1, t}$. Without knowing whether the change occurs and where it occurs beforehand when forming the GLR statistic, we have to maximize $k$ over all possible change locations. The GLR statistic is given by $\max_{k<t} \sum_{i=k+1}^t \log(f_{\hat{\theta}_{k, t}}(X_i)/f_{\theta_0}(X_t))$, and a change is announced whenever it exceeds a pre-specified threshold. The GLR statistic is more robust than CUSUM \citep{lai1998information}, and it is particularly useful when the post-change parameter may vary from one situation to another. \yao{In simple cases, the MLE $\hat{\theta}_{k+1, t}$ may have closed-form expressions and may be evaluated recursively. For instance, when the post-change distribution is Gaussian with mean $\theta$ \cite{lorden2005nonanticipating}, $\hat{\theta}_{k+1, t} = (\sum_{i=k+1}^t X_{i})/(t-k)$, and $\hat{\theta}_{k+1, t+1} = (t-k)/(t-k+1)\cdot\hat{\theta}_{k+1, t} + X_{t+1}/(t-k+1)$.} However, \yao{in more complex situations, in general MLE $\hat{\theta}_{k+1, t}$ does not have recursive form and cannot be evaluated using simple summary statistics. One such instance is given in Section \ref{sec:social_net}. Another instance is when there is a constraint on the MLE such as sparsity.} In these cases, one has to store historical data and \yc{recompute the MLE} $\hat{\theta}_{k, t}$  whenever there is new data, which is not memory efficient nor computational efficient. For these cases, as a remedy, the window-limited GLR is usually considered, where only the past $w$ samples are stored and the maximization is restricted to be over $k\in(t-w, t]$. However, even with the window-limited GLR, one still has to recompute $\hat{\theta}_{k, t}$ using historical data whenever the new data are added. 

Besides CUSUM or GLR, various online change-point detection procedures using one-sample updates have been considered, which replace with the MLE with a simple recursive estimator.  The one-sample update estimate takes the form of $\hat{\theta}_{k, t} =  h(X_t, \hat{\theta}_{k, t-1})$ for some function $h$ that uses only the most recent data and the previous estimate. Then the estimates are plugged into the likelihood ratio statistic to perform detection. \yao{Online convex optimization algorithms (such as online mirror descent) are natural approach to construct these estimators (see, e.g.,  \cite{raginsky2009sequential,raginsky2012sequential}). Such a scheme provides a more versatile approach to develop detecting procedure for complex situations, where the exact MLE does not have a recursive form or even a closed-form expression.} The one-sample update enjoys efficient computation, as information from the new data can be incorporated via low computational cost update. 
It is also memory efficient since the update only needs the most recent sample. 
The one sample update estimators may not correspond to the exact MLE, but they tend to result in good detection performance. However, in general there is no performance guarantees for such approach. This is the question we aim to address in this paper.


\subsection{Application scenario: Social network change-point detection}\label{sec:social_net}

The widespread use of social networks (such as Twitter) leads to a large amount of user-generated data generated continuously. One important aspect is to detect change-points in streaming social network data. These change-points may represent the collective anticipation of  response to external events or system ``shocks'' \cite{PeelClauset2014}. Detecting such changes can provide a better understanding of patterns of social life. In social networks, a common form of the data is discrete events over continuous time. As a simplification, each event contains a time label and a user label in the network. In our prior work \cite{LiXie17}, we model discrete events  using network point processes, which capture the influence between users through an {\it influence matrix}. We then cast the problem as detecting changes in \yc{an} influence matrix, assuming that the influence matrix in the normal state (before the change) can be estimated from the reference data. After the change, the influence matrix is unknown (since it represents an anomaly) and has to be estimated online. Due to computational burden and memory constraint, since the scale of the network tends to be large, we do not want to store the entire historical data and rather compute the statistic in real-time. A simulated example to illustrate this case is shown later in Section \ref{sec:poisson_hawkes}. 

\subsection{Contributions}

\yao{This paper has two main contributions. First, we present a general approach based on online convex optimization (OCO) for constructing the estimator for the one-sided sequential hypothesis test and the sequential change-point detection, in the non-anticipative approach of \cite{lorden2005nonanticipating} if the MLE cannot be computed in a convenient recursive form.}

\yao{Second, we provide a proof of the near second-order asymptotic optimality of this approach when a ``logarithmic regret property'' is satisfied and when the distributions are from an exponential family.} 
The nearly second-order asymptotic optimality \citep{tartakovsky2014sequential} means that the upper bound for performance matches the lower bound up to a log-log factor as the false-alarm rate tends to zero. Inspired by the existing connection between sequential analysis and online convex optimization in \cite{cesa2006prediction, hazan2016introduction}, we prove the near optimality leveraging the logarithmic regret property of online mirror descent (OMD) and the lower bound established in statistical sequential change-point literature \citep{siegmund2008minimax,tartakovsky2014sequential}. More precisely, we provide a general upper bound for  one-sided sequential hypothesis test and change-point detection procedures with the one-sample update schemes. The upper bound explicitly captures the impact of estimation on detection by an  estimation algorithm dependent factor. This factor shows up as an additional term in the upper bound for the expected detection delay, and it corresponds to the regret incurred by the one-sample update estimators. \footnote{This establishes an interesting linkage between sequential change-point detection and online convex optimization. Although both fields, sequential change-point detection and online convex optimization, study sequential data, the precise connection between them is not clear, partly because the performance metrics are different: the former concerns with the tradeoff between average run length and detection delay, whereas the latter focuses on bounding the cumulative loss incurred by the sequence of estimators through a regret bound \citep{azoury2001relative,hazan2016introduction}.} Synthetic examples validate the performances of one sample update schemes. Here we focus on OMD estimators, but the results can be generalized to other OCO schemes such as the online gradient descent.

\subsection{Literature and related work}

Sequential \yc{change-point} detection is a classic subject with an extensive literature. Much success has been achieved when the pre-change and post-change distributions are exactly specified. For example, the CUSUM procedure \citep{page1954continuous} with first-order asymptotic optimality \cite{lorden1971procedures} and exact optimality \cite{moustakides1986optimal} in the minimax sense, \yc{and }the Shiryayev-Roberts (SR) procedure \cite{shiryaev1963optimum} derived based on  Bayesian principle \yc{that also} enjoys various optimality \yc{properties}. Both CUSUM and SR procedures rely on likelihood ratios between the specified pre-change and post-change distributions. 

\yao{There are two main approaches in dealing with the unknown post-change parameters. The first one is a GLR approach \cite{willsky1976generalized,lai1995sequential,lai1998information,lai2004likelihood,lorden2008sequential}, and the second is a mixture approach \cite{siegmund2008minimax,pollak1987average}. The GLR statistic enjoys certain optimality properties, but it can not be computed recursively in many cases \cite{lai2004likelihood}. To address the infinite memory issue, \cite{willsky1976generalized,lai1998information} studied the window-limited GLR procedure. The main advantage of the mixture approach is that it allows an easy evaluation of a threshold that guarantees the desired false alarm constraint. A disadvantage of this approach is that sometimes there may not be a natural way of selecting the weight function, in particular when there is no conjugate prior. This motivated a third approach to this problem, which was proposed first by Robbins and Siegmund in the context of hypothesis testing, and then Lorden and Pollak \cite{lorden2005nonanticipating} in the sequential change detection problem. This approach replaces the unknown parameter with some non-anticipating estimator, which can be easier to find even if there is no conjugate prior, as in the Gamma example considered in \cite{lorden2005nonanticipating,pollak1987average}. These work developed a modified SR procedure by introducing a prior distribution to the unknown parameters. 
}
%
%
\yao{While the non-anticipating estimator approach \cite{lorden2005nonanticipating,lorden2008sequential} enjoys recursive and thus efficient computation for the likelihood ratio based detection statistics, but their approaches to construct recursive estimators (based on MLE or method-of-moments) cannot be easily extended to more complex cases  (for instance, multi-dimensional parameters with constraints).
Here, we consider a general and convenient approach for constructing non-anticipating estimators based on online convex optimization which is particularly useful for these complex cases}. 
Our work provides an alternative proof for the nearly second-order asymptotic optimality by building a connection to online convex optimization and leveraging the regret bound type of results  \cite{hazan2016introduction}. For one-dimensional Gaussian mean shift without any constraint, we replicate the second-order asymptotic optimality, namely, Theorem 3.3 in \cite{lorden2008sequential}. \yao{Recent work \cite{Mei06} also treats the problem when the pre-change distribution has unknown parameters.}

Another related problem is sequential joint estimation and detection, but the goal is different in that one aims to achieve both good detection and good estimation performance, whereas in our setting estimation is only needed for computing the detection statistics. These works include 
\cite{yilmaz2015sequential} and \cite{yilmaz2016sequential}, which study the joint detection and estimation problem of a specific form that arises from many applications such as spectrum sensing \cite{yilmaz2014sequential}, image observations \cite{vo2010joint}, and MIMO radar \cite{tajer2010optimal}: a linear scalar observation model with Gaussian noise, and under the alternative hypothesis there is {\it an unknown multiplicative parameter}. \yc{The paper} of \cite{yilmaz2015sequential} demonstrates that solving the joint problem by treating detection and estimation separately with the corresponding optimal procedure does not yield an overall optimum performance, and provides an elegant closed-form optimal detector. Later on \cite{yilmaz2016sequential} generalizes the results. There are also other approaches solving the joint detection-estimation problem using multiple hypotheses testing \citep{baygun1995optimal,vo2010joint} and Bayesian formulation\yc{s} \cite{moustakides2012joint}.

Related work using online convex optimization for anomaly detection includes \cite{raginsky2009sequential}, which develops an efficient detector for the exponential family using online mirror descent and proves a logarithmic regret bound, and \cite{raginsky2012sequential}, which dynamically adjusts the detection threshold to allow feedbacks about whether decision outcome. However, these works consider a different setting that the change is a transient outlier instead of a persistent change, as assumed by the classic statistical change-point detection literature. When there is persistent change, it is important to accumulate ``evidence'' by pooling the post-change samples (our work considers the persistent change).

Extensive work has been done for parameter estimation in the online-setting. This includes online density estimation over the exponential family by regret minimization \citep{azoury2001relative,raginsky2009sequential,raginsky2012sequential}, sequential prediction of individual sequence with the logarithm loss \citep{cesa2006prediction,kotlowski2011maximum}, online prediction for time series \cite{Hazan-time-series-13}, and sequential NML (SNML) prediction \cite{kotlowski2011maximum} which achieves the optimal regret bound. Our problem is different from the above, in that estimation is not the end goal; one only performs parameter estimation to plug them back into the likelihood function for detection. Moreover, a subtle but important difference of our work is that the loss function for online detecting estimation is $-f_{\hat{\theta}_i}(X_i)$, whereas our loss function is $-f_{\hat{\theta}_{i-1}}(X_i)$ in order to retain the {\it martingale property}, which is essential to establish the nearly second-order \yc{asymptotic} optimality.


\section{Preliminaries}

Assume a sequence of i.i.d. random variables $X_1, X_2, \ldots$ with a probability density function of a parametric form $f_{\theta}$. The parameter $\theta$ may be unknown. Consider two related problems: \yc{one-sided} sequential hypothesis test and sequential change-point detection. The detection statistic relies on a sequence estimators $\{\hat{\theta}_t\}$ constructed using online mirror descent. The OMD uses simple {\it one-sample update}: the update from  $\hat{\theta}_{t-1}$ to $\hat{\theta}_{t}$ only uses the current sample $X_t$. This is the main difference from the traditional generalized likelihood ratio (GLR) statistic \cite{lai1998information}, where each  $\hat{\theta}_t$ is estimated using historical samples. In the following, we present detailed descriptions for two problems. We will consider exponential family \yc{distributions} and present our non-anticipating estimator based on the one-sample estimate.

\subsection{\yc{One-sided} sequential hypothesis test}
\yc{First, we consider a one-sided sequential hypothesis test where the goal is only to reject the null hypothesis. This is a special case of the change-detection problem where the change-point can be either $0$ or $\infty$ (meaning it never occurs). Studying this special case will given us an important intermediate step towards solving the sequential change-detection problem.} 

Consider \yc{the }null hypothesis $\textsf{H}_0: \theta = \theta_0$ versus the alternative $\textsf{H}_1: \theta \neq \theta_0$. Hence the parameter under the alternative distribution is unknown. The classic approach to solve this problem is the \yc{one-sided }sequential probablity-ratio test (SPRT) \citep{wald1948optimum}: at each time, given samples $\{X_1, X_2, \ldots, X_{t}\}$, the decision is either to \yc{reject $\textsf{H}_0$} or taking more samples if \yc{the rejection decision} can\yc{not }be \yc{made} confidently. Here, we introduce \yc{a }{\it modified} \yc{one-sided }SPRT with a sequence of {\it non-anticipating} plug-in estimators:
\begin{equation}
\hat{\theta}_t := \hat{\theta}_t(X_1, \ldots, X_{t}), \quad t = 1, 2, \ldots.
\label{eq:seq_estimator}
\end{equation} 
Define the \yc{test statistic} at time $t$ as
\begin{equation}
\Lambda_t = \prod_{i=1}^t \frac{f_{\hat{\theta}_{i-1}}(X_i)}{f_{\theta_0}(X_i)}, \quad  i\geq 1. 
\label{onestat}
\end{equation} 
\yc{The} test statistic has a simple recursive implementation: 
\[\Lambda_t = \Lambda_{t-1} \cdot \frac{f_{\hat{\theta}_{t-1}}(\yc{X_t})}{f_{\theta_0}(\yc{X_t}}.\] \yc{Define a sequence of $\sigma$-algebras $\{\mathcal{F}_t\}_{t\geq 1}$ where $\mathcal{F}_t = \sigma(X_1, \ldots, X_t)$}. \yc{The test statistic has the martingale property due to its non-anticipating nature: $\mathbb{E}[\Lambda_t \mid \mathcal{F}_{t-1}] =  \Lambda_{t-1}$, where the expectation is taken when $X_1,\ldots$ are i.i.d. random variables drawn from $f_{\theta_0}$.} The decision rule is a stopping time 
\begin{equation}
\tau(b) = \min\{t\geq 1: \log \Lambda_t \geq b\},
\label{sequentialHP}
\end{equation}
where $b>0$ is a pre-specified threshold. We reject the null hypothesis whenever the statistic exceeds the threshold. The goal is to \yc{reject the null hypothesis} using as few samples as possible under the \yc{false-alarm rate (or Type-I error)} constraint.

\subsection{Sequential change-point detection} 
\yc{Now we consider the sequential change-point detection problem.} A change may occur at an unknown time $\nu$ which alters the underlying distribution of the data. One would like to detect such a change as quickly as possible. Formally, change-point detection can be cast into the following hypothesis test: 
\begin{equation}
\begin{split}
\textsf{H}_0: &~~X_1, X_2, \ldots \overset{\rm i.i.d.}{\sim} f_{\theta_0}, \\ 
\textsf{H}_1: &~~X_1, \ldots, X_{\nu} \overset{\rm i.i.d.}{\sim} f_{\theta_0}, \quad X_{\nu+1}, X_{\nu+2}, \ldots \overset{\rm i.i.d.}{\sim} f_{\theta}, 
\end{split}
\label{maintestproblem}
\end{equation}
Here we assume \yc{an unknown }$\theta$ \yc{to represent the anomaly.} The goal is to detect the change as quickly as possible after it occurs under the false-alarm rate constraint. 
We will consider likelihood ratio based detection procedures adapted from two types of existing ones, which we call \yc{the} adaptive CUSUM (ACM), and the adaptive SRRS (ASR) procedures.

For change-point detection, the post-change parameter is estimated using post-change samples. This means that, for each putative change-point location before the current time $k<t$, the post-change samples are $\{X_{k}, \ldots, X_t\}$; with a slight abuse of notation, the post-change parameter is estimated as 
\begin{equation}
\hat{\theta}_{k,i} = \hat{\theta}_{k,i}(X_k, \ldots, X_{i}), \quad i\geq k.
\label{eq:theta_change}
\end{equation} 
Therefore, for $k=1$, $\hat{\theta}_{k,i}$ becomes $\hat{\theta}_i$ defined in (\ref{onestat}) for \yc{the one-sided }SPRT. Initialize with $\hat{\theta}_{k,k-1} = \theta_0$. The likelihood ratio at time $t$ for a hypothetical change-point location $k$ is given by 
\begin{equation}
\Lambda_{k,t} = \prod_{i=k}^t  \frac{f_{\hat{\theta}_{k,i-1}}(X_i)}{f_{\theta_0}(X_i)}, 
\label{cumulativestat}
\end{equation}
where $\Lambda_{k, t}$ can be computed recursively similar to (\ref{onestat}).

\begin{figure}[h]
\begin{center}
\includegraphics[width=0.7\textwidth]{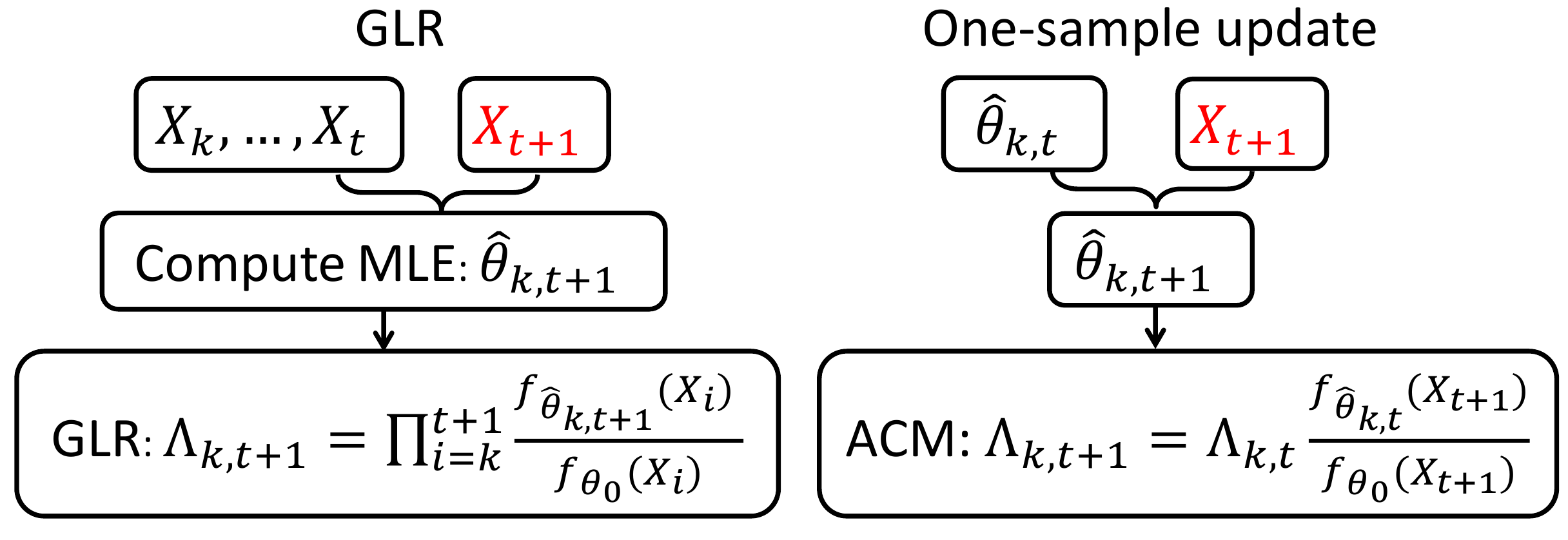}
\end{center}
\vspace{-0.1in}
\caption{Comparison of the update scheme for GLR and our methods when a new sample arrives.}
\label{fig:comparison_GLRACM}
\end{figure}

Since we do not know the change-point location $\nu$, from the maximum likelihood principle, we take the maximum of the statistics over all possible values of $k$. This gives the ACM procedure:
\begin{equation}
T_{\rm ACM}(\yc{b_1}) = \inf \left\{t\geq 1: \max_{1\leq k\leq t} \log \Lambda_{k,t} > \yc{b_1}\right\}, 
\label{ACMprocedure}
\end{equation}
where $\yc{b_1}$ is a pre-specified threshold. 
Similarly, by replacing the maximization \yc{over $k$} in (\ref{ACMprocedure}) with summation, we obtain the following ASR procedure \citep{lorden2005nonanticipating}, which can be interpreted as a Bayesian statistic similar to the Shiryaev-Roberts procedure. 
\begin{equation}
T_{\rm ASR}(\yc{b_2}) = \inf \left\{t\geq 1: \log \left(\sum_{k=1}^t \Lambda_{k,t} \right)> \yc{b_2}  \right\},
\label{ASRprocedure}
\end{equation}
where $\yc{b_2}$ is a pre-specified threshold. The computations of $\Lambda_{k,t}$ and estimator $\{\hat{\theta}_t\}$, $\{\hat{\theta}_{k,t} \}$ are discussed later in section \ref{sec:estimators}. For a fixed $k$, the comparison between our methods and GLR is illustrated in Figure \ref{fig:comparison_GLRACM}.

\vspace{-0.1in}
\begin{Remark}
In practice, to prevent the memory and computation complexity from blowing up as time $t$ goes to infinity, we can use window-limited version of the detection procedures in (\ref{ACMprocedure}) and (\ref{ASRprocedure}). The window-limited versions are obtained by replacing $\max_{1\leq k\leq t}$ with $\max_{t-w\leq k\leq t}$ in (\ref{ACMprocedure}) and by replacing $\sum_{k=1}^t$ with $\sum_{k=t-w}^t$ in (\ref{ASRprocedure}). Here $w$ is a prescribed window size. Even if we do not provide theoretical analysis to the window-limited versions, we refer the readers to \cite{lai1998information} for the choice of $w$ the window-limited GLR procedures. 
\label{remark_window}
\end{Remark}

\subsection{Exponential family}
\label{sec:ef}

In this paper, we focus on $f_\theta$ being the exponential family for the following reasons: (i) exponential family \citep{raginsky2012sequential} represents a very rich class of parametric and even many nonparametric statistical models \citep{barron1991approximation}; (ii) the negative log-likelihood function for exponential family $-\log f_{\theta}(x)$ is convex, and this allows us to perform online convex optimization. Some useful properties of the exponential family are briefly summarized below, and full proofs can be found in \cite{wainwright2008graphical, raginsky2012sequential}.

Consider an observation space $\mathcal{X}$ equipped with a sigma algebra $\mathcal{B}$ and a sigma finite measure $H$ on $(\mathcal{X},\mathcal{B})$. Assume the number of parameters is $d$. Let $x^\intercal$ denote the transpose of a vector or matrix. Let $\phi: \mathcal{X}\rightarrow \mathbb{R}^d$ be an $H$-measurable function $\phi(x) = (\phi_1(x),\ldots, \phi_d(x))^\intercal$. Here $\phi(x)$ corresponds to the  sufficient statistic for $\theta$. Let $\Theta$ denote the parameter space in $\mathbb{R}^d$. Let $\{\mathcal{P}_{\theta}, \theta \in \Theta \}$ be a set of probability distributions with respect to the measure $H$. Then, $\{\mathcal{P}_{\theta}, \theta \in \Theta \}$ is said to be a multivariate exponential family with natural parameter $\theta$, if the probability density function of each $f_{\theta} \in \mathcal{P}_{\theta}$ with respect to $H$ can be expressed as
$
f_{\theta}(x) = \exp\{ \theta^\intercal \phi(x) - \Phi(\theta) \}.  
$
In the definition, the so-called log-partition function is given by
\[
\Phi(\theta) := \log \int_{\mathcal{X}} \exp(\theta^\intercal \phi(x)) d H(x).
\]
To make sure $f_{\theta}(x)$ a well-defined probability density, we consider the following two sets for parameters:
\[
\Theta = \{ \theta \in \mathbb{R}^d: \log \int_{\mathcal{X}} \exp(\theta^\intercal \phi(x)) d H(x) < +\infty \},\]
and
\[\Theta_\sigma = \{\theta \in \Theta: \nabla^2 \Phi(\theta) \succeq \sigma I_{d\times d} \}.
\]
Note that $-\log f_{\theta}(x)$ is $\sigma$-strongly convex over $\Theta_\sigma$.
Its gradient corresponds to
$
\nabla\Phi(\theta) = \mathbb{E}_{\theta}[\phi(X)]
$, and the Hessian $\nabla^2\Phi(\theta)$ corresponds to the covariance matrix of the vector $\phi(X)$. \yc{Therefore,} \yc{$\nabla^2$}$\Phi(\theta)$ is positive semidefinite and $\Phi(\theta)$ is convex. 
Moreover, $\Phi$ is a {\it Legendre function}, which means that it is strongly convex, continuous differentiable and essentially smooth \cite{wainwright2008graphical}. 
The Legendre-Fenchel dual $\Phi^*$ is defined as 
\[
\Phi^*(z) = \sup_{u \in \Theta}\{u^\intercal z - \Phi(u) \}.
\] The mappings $\nabla \Phi^*$ is an inverse mapping of  $\nabla \Phi$ \citep{beck2003mirror}. Moreover, if $\Phi$ is a strongly convex function, then $\nabla \Phi^* = (\nabla \Phi)^{-1}$. 

A general measure of proximity used in \yc{the }OMD is the so-called {\it Bregman divergence} $B_F$, which is a nonnegative function induced by a Legendre function $F$ (see, e.g., \citep{wainwright2008graphical, raginsky2012sequential}) defined as
\begin{equation}
B_F(u, v):= F(u) - F(v) - \langle \nabla F(v), u-v \rangle.\label{eq:Bregman}
\end{equation}
For exponential family, a natural choice of the Bregman divergence is the Kullback-Leibler (KL) divergence. Define $\mathbb{E}_{\theta}$ as the expectation when $X$ is a random variable with density $f_{\theta}$ and $I(\theta_1, \theta_2)$ as the KL divergence between two distributions with densities $f_{\theta_1}$ and $f_{\theta_2}$ for any $\theta_1, \theta_2 \in \Theta$. Then
\begin{equation}
I(\theta_1, \theta_2) = \mathbb{E}_{\theta_1}\left[\log (f_{\theta_1}(X)/f_{\theta_2}(X)) \right]. 
\label{KL}
\end{equation}
It can be shown that, for exponential family, 
$
I(\theta_1, \theta_2) = \Phi(\theta_2) - \Phi(\theta_1) - (\theta_2-\theta_1)^\intercal \nabla\Phi(\theta_1).
$
Using the definition (\ref{eq:Bregman}), this means that $B_{\Phi}$
\begin{equation}
B_{\Phi}(\theta_1, \theta_2) := I(\theta_2, \theta_1)
\label{eq:bregmean_exponential}
\end{equation}
is a Bregman divergence. This property is  useful to constructing mirror descent estimator for the exponential family   \citep{nemirovskii1983problem,beck2003mirror}.

\subsection{Online convex optimization (OCO) algorithms for non-anticipating estimators} 
\label{sec:estimators}

\ycc{Online convex optimization (OCO) algorithms \cite{hazan2016introduction} can be interpreted as a player who  makes sequential decisions. At the time of each decision, the  outcomes are unknown to the player.
After committing to a decision, the decision maker suffers a loss that can be adversarially chosen. An OCO algorithm makes decisions, which, based on the observed outcomes, minimizes the
{\it regret} that  is the difference between the total loss that has incurred relatively to that of the best fixed decision in hindsight. }
\ycc{To design non-anticipating estimators, we consider OCO algorithms with  likelihood-based regret functions. We iteratively estimate the parameters at the time when a one new observation becomes available based on {\it the maximum likelihood principle}, and hence the loss incurred corresponds to the negative log-likelihood of the new sample evaluated at the estimator} $\ell_t(\theta) := -\log f_{\theta}(X_t)$, \yc{which corresponds to the log-loss in \cite{cesa2006prediction}}. 
Given samples $X_1, \ldots, X_t$, the regret for a sequence of estimators $\{\hat{\theta}_i\}_{i=1}^t$ generated by a {\it likelihood-based OCO algorithm} $\textsf{a}$ is defined as
\begin{equation}
\mathcal R_t^{\textsf{a}} = \sum_{i=1}^t \{-\log f_{\hat{\theta}_{i-1}}(X_i)\} - \inf_{\tilde{\theta} \in \Theta} \sum_{i=1}^t \{-\log f_{\tilde{\theta}}(X_i)\}.
\label{regret}
\end{equation}
Below we omit the superscript $\textsf{a}$ occasionally for notational simplicity.

In this paper, we consider a generic OCO procedure called the online mirror descent algorithms (OMD) \cite{hazan2016introduction,shalev2012online}. 
 Next, we discuss how to construct the non-anticipating estimators $\{\hat{\theta}_t\}_{t\geq 1}$ in (\ref{eq:seq_estimator}), and $\{\hat{\theta}_{k, t}\}, k = 1, 2, \ldots, t-1$ in (\ref{eq:theta_change}) using \yc{OMD}. 
The main idea of OMD is the following. At each time step, the estimator $\hat{\theta}_{t-1}$ is updated using the new sample $X_t$, by balancing the tendency to stay close to the previous estimate against the tendency to move in the direction of the greatest local decrease of the loss function. For the loss function defined above, a sequence of OMD estimator is constructed by
\begin{equation}
\hat{\theta}_{t} = \mathop{\arg\min}_{u \in \Gamma} [ u^\intercal \nabla \ell_t(\hat{\theta}_{t-1}) + \frac{1}{\eta_i}B_{\Phi}(u, \hat{\theta}_{t-1}) ],
 \label{Breg2}
\end{equation} 
\yc{where $B_{\Phi}$ is defined in (\ref{eq:bregmean_exponential}).} 
Here $\Gamma \subset \Theta_\sigma$ is a closed convex set, which is problem-specific and encourages certain parameter structure such as sparsity. 
\begin{Remark}
\label{remark:CD}
Similar to (\ref{Breg2}), for any fixed $k$, we can compute $\{\hat{\theta}_{k, t}\}_{t\geq 1}$ via OMD for sequential change-point detection. The only difference is that $\{\hat{\theta}_{k, t}\}_{t\geq 1}$ is computed if we use $X_k$ as our first sample and then apply the recursive update (\ref{Breg2}) on $X_{k+1}, \ldots$. For $\hat{\theta}_{t}$, we use $X_1$ as our first sample.
\end{Remark}

There is an equivalent form of OMD, presented as the original formulation in \cite{nemirovskii1983problem}. The equivalent form is sometimes easier to use for algorithm development, and it consists of four steps: (1) compute the dual variable: $\hat{\mu}_{t-1} = \nabla \Phi(\hat{\theta}_{t-1})$; (2) perform the dual update: $\hat{\mu}_{t} = \hat{\mu}_{t-1} - \eta_t\nabla \ell_t(\hat{\theta}_{t-1})$; (3) compute the primal variable: $\tilde{\theta}_t = (\nabla \Phi)^{*}(\hat{\mu}_t)$; (4) perform the projected primal update: $\hat{\theta}_t = \mathop{\arg\min}_{u \in \Gamma} B_{\Phi}(u, \tilde{\theta}_t)$. The equivalence between the above form for OMD and the nonlinear projected subgradient approach in (\ref{Breg2}) is proved in \cite{beck2003mirror}. We adopt this approach when deriving our algorithm and follow the same strategy as \cite{raginsky2009sequential}. Algorithm \ref{alg1} summarizes the steps\footnote{The implementation of the code can be downloaded at \url{http://www2.isye.gatech.edu/~yxie77/one-sample-update-code.zip}.}.

\begin{algorithm}[h!]
\caption{Online mirror-descent for \yc{non-anticipating} estimators
}\label{alg1}
\begin{algorithmic}[1]

\REQUIRE Exponential family specifications $\phi(x), \Phi(x)$ and $f_{\theta}(x)$; initial parameter value $\theta_0$; sequence of data $X_1, \ldots, X_t, \ldots$; a closed, convex set for parameter $\Gamma \subset\Theta_\sigma$; a decreasing sequence $\{ \eta_t\}_{t\geq 1}$ of strictly positive step-sizes.

\vspace{.1in}
\STATE $\hat{\theta}_0 = \theta_0, \Lambda_0 = 1$. \COMMENT{Initialization}

\vspace{.1in}
\FORALL{$t = 1,2,\ldots,$}

\STATE \mbox{Acquire a new observation $X_t$}

\vspace{.1in}
\STATE \mbox{Compute loss $\ell_t(\hat{\theta}_{t-1}) \triangleq -\log f_{\hat{\theta}_{t-1}}(X_t) = \Phi(\hat{\theta}_{t-1}) - \hat{\theta}_{t-1}^\intercal \phi(X_t)$}

\vspace{.1in}
\STATE 
\mbox{Compute likelihood ratio} $
\Lambda_t = \Lambda_{t-1} \dot f_{\hat{\theta}_{t-1}}(X_t)/f_{\theta_0}(X_t)
$ 


\vspace{.1in}
\STATE $\hat{\mu}_{t-1} = \nabla \Phi(\hat{\theta}_{t-1})$,  $\hat{\mu}_{t} 
 = \hat{\mu}_{t-1} - \eta_t(\hat{\mu}_{t-1} - \phi(X_{t}))$ \COMMENT{Dual update}

\vspace{.1in}
\STATE $\tilde{\theta}_t = (\nabla \Phi)^{*}(\hat{\mu}_t)$

\vspace{.1in}
\STATE  $\hat{\theta}_t = \mathop{\arg\min}_{u \in \Gamma} B_{\Phi}(u, \tilde{\theta}_t)$  \COMMENT{Projected primal update}
\vspace{.1in}
\ENDFOR
\vspace{.1in}
\RETURN $\{\hat{\theta}_t\}_{t\geq 1}$ and $\{\Lambda_t\}_{t\geq 1}$.

\end{algorithmic}

\end{algorithm}

For strongly convex loss function, the regret of many OCO algorithms, including the \yc{OMD}, has the property that $\mathcal{R}_n \leq C \log n $ for some constant $C$ (depend on $f_\theta$ and $\Theta_\sigma$) and any positive integer $n$ \citep{AgarwalDuchi,raginsky2012sequential}. Note that for exponential family, the loss function is the negative log-likelihood function, which is strongly convex over $\Theta_{\sigma}$. 
Hence, we can have the logarithmic regret property.

\section{Nearly second-order \yc{asymptotic} optimality of one-sample update schemes}
\label{sec:analysis}

Below we prove the {\it nearly second-order \yc{asymptotic} optimality} of the one-sample update schemes. More precisely, the nearly second-order \yc{asymptotic }optimality means that the algorithm obtains the lower performance bound asymptotically up to a log-log factor in the false-alarm rate, as the false-alarm rate tends to zero (in many cases the log-log factor is a small number).

We first introduce some necessary notations. Denote $\mathbb{P}_{\theta,\nu}$ and $\mathbb{E}_{\theta, \nu}$ \yc{as} the probability measure and \yc{the} expectation when the change occurs at time $\nu$ and the post-change parameter is $\theta$, i.e., when $X_1,\ldots, X_{\nu}$ are i.i.d. random variables with density $f_{\theta_0}$ and $X_{\nu+1},X_{\nu+2},\ldots$ are i.i.d. random variables with density $f_{\theta}$. Moreover, let $\mathbb{P}_{\infty}$ and $\mathbb{E}_{\infty}$ denote the probability measure when there is no change, i.e., $X_1, X_2, \ldots$ are i.i.d. random variables with density $f_{\theta_0}$. Finally, let $\mathcal{F}_t$ denote the $\sigma$-algebra generated by $X_1, \ldots, X_t$ for $t\geq 1$.

\subsection{\yc{``One-sided''} Sequential hypothesis test}
\label{sec:onesided}

\yao{Recall that the decision rule for sequential hypothesis test is a stopping time $\tau(b)$ defined in (\ref{sequentialHP}).}
The two standard performance metrics are the \yc{false-alarm rate}, \yc{denoted} as $\mathbb{P}_{\infty}(\tau(b) < \infty)$, \yc{and the expected detection delay (i.e., the expected number of samples needed to reject the null)}\yc{, denoted as} $\mathbb{E}_{\theta, 0}[\tau(b)]$. A meaningful test should have both small $\mathbb{P}_{\infty}(\tau(b) < \infty)$ and small $\mathbb{E}_{\theta, 0}[\tau(b)]$. Usually, one adjusts the threshold $b$ to control the \yc{false-alarm rate} to be below a certain level. 

Our main result is the following. As has been observed by \cite{lai2004likelihood}, there is a loss in the statistical efficiency by using one-sample update estimators relative to the GLR approach using the entire sample\yc{s} $X_1, \ldots, X_t$ in the past. The theorem below shows that this loss corresponds to the expected regret \yc{given in (\ref{regret}).}

\begin{Theorem}[Upper bound for \ycc{OCO} based SPRT]
\ycc{Let $\{\hat{\theta}_t \}_{t\geq 1}$ be a sequence of non-anticipating estimators generated by an OCO algorithm $\textsf{a}$. } As $b\rightarrow \infty$,
\begin{equation}
\mathbb{E}_{\theta, 0}[\tau(b)] \leq \frac{b}{I(\theta, \theta_0)} + \frac{ \mathbb{E}_{\theta,0}\left[ \mathcal R^\textsf{a}_{\tau(b)}\right]}{I(\theta, \theta_0)}+O(1)
\label{sequentialHPdelay}
\end{equation}
Here $O(1)$ is a term upper-bounded by an absolute constant as $b\rightarrow \infty$.
\label{maintheorem}
\end{Theorem}

The main idea of the proof is to decompose the statistic defining $\tau(b)$, $\log\Lambda(t)$, into a few terms that form martingales, and then \yc{invoke} the Wald's Theorem for the stopped process. 
\yc{
\begin{Remark} The inequality (\ref{sequentialHPdelay}) is valid for a sequence of non-anticipating estimators generated by an OCO algorithm. Moreover, (\ref{sequentialHPdelay}) gives an explicit connection between the expected detection delay for the one-sided sequential hypothesis testing (left-hand side of (\ref{sequentialHPdelay})) and the regret for the OCO (the second term on the right-hand side of (\ref{sequentialHPdelay})). This illustrates clearly the impact of estimation on detection by an estimation algorithm dependent factor.
\end{Remark}
}

Note that in the statement of the Theorem \ref{maintheorem}, the stopping time $\tau(b)$ appears on the right-hand side of the inequality (\ref{sequentialHPdelay}). For OMD, the expected sample size is usually small. By comparing with specific regret bound $\mathcal R_{\tau(b)}$, we can bound $\mathbb{E}_{\theta,0}[\tau(b)]$ as discussed in Section \ref{sec:simulation}. The most important case is that when the estimation algorithm has a logarithmic expected regret. For the exponential family, as shown in section \ref{sec:regretboundexample}, Algorithm \ref{alg1} can achieve $\mathbb{E}_{\theta,0}[\mathcal R_n] \leq C\log n$ for any positive integer $n$. \yc{To obtain a more specific order of the upper bound for $\mathbb{E}_{\theta,0}[\tau_b]$ when  $b$ grows, we establish an upper bound for $\mathbb{E}_{\theta,0}[\tau_b]$ as a function of $b$,} to obtain the following Corollary \ref{cor:logregret}.

\begin{Corollary}
\ycc{Let $\{\hat{\theta}_t \}_{t\geq 1}$ be a sequence of non-anticipating estimators generated by an OCO algorithm $\textsf{a}$.} Assume that $\mathbb{E}_{\theta,0}[\mathcal R_{n}^\textsf{a}] \leq C\log n$ for any positive integer $n$ and some constant $C>0$, we have
\begin{equation}
\mathbb{E}_{\theta,0}[\tau(b)] \leq \frac{b}{I(\theta, \theta_0)} + \frac{C\log b}{I(\theta, \theta_0)}(1+o(1)).
\label{logregretresult}
\end{equation}
Here $o(1)$ is a vanishing term as $b \rightarrow \infty$. 
\label{cor:logregret}
\end{Corollary}
\yc{Corollary \ref{cor:logregret} shows that other than the well known first-order approximation $b/I(\theta, \theta_0)$ \cite{lorden1971procedures, lorden2005nonanticipating}, the expected detection delay $\mathbb{E}_{\theta,0}[\tau(b)]$ is bounded by an additional term that is on the order of $\log(b)$ if the estimation algorithm has a logarithmic regret. This $\log b$ term plays an important role in establishing the optimality properties later. 
To show the optimality properties for the detection procedures, we first select a set of detection procedures with false-alarm rates lower than a prescribed value, and then prove that among all the procedures in the set, the expected detection delays of our proposed procedures are the smallest. 
Thus, we can choose a threshold $b$ to uniformly control the false-alarm rate of $\tau(b)$.} 
\begin{Lemma}[false-alarm rate of $\tau(b)$]
\yc{Let $\{\hat{\theta}_t \}_{t\geq 1}$ be any sequence of non-anticipating estimators. For any $b>0$, } $\mathbb{P}_{\infty}(\tau(b) < \infty) \leq \exp(-b)$.
\label{typeIerrorbound}
\end{Lemma}

Lemma \ref{typeIerrorbound} \yc{shows that as $b$ increases the false-alarm rate of $\tau(b)$ decays exponentially fast. We can set $b=\log(1/\alpha)$ to make the false-alarm rate of $\tau(b)$  less than some $\alpha>0$.} Next, leveraging an existing lower bound for general SPRT presented in Section 5.5.1.1 in \citep{tartakovsky2014sequential}, we establish the nearly second-order \yc{asymptotic }optimality of OMD based SPRT as follows: 

\begin{Corollary}[Nearly second-order optimality of \ycc{OCO} based SPRT] \ycc{Let $\{\hat{\theta}_t \}_{t\geq 1}$ be a sequence of non-anticipating estimators generated by an OCO algorithm $\textsf{a}$.} Assume that $\mathbb{E}_{\theta,0}[\mathcal R^\textsf{a}_{n}] \leq C\log n$ for any positive integer $n$ and some constant $C>0$. Define a set $C(\alpha) = \{T: \mathbb{P}_{\infty}(T < \infty) \leq \alpha \}$. For $b=\log(1/\alpha)$, due to Lemma \ref{typeIerrorbound}, $\tau(b) \in C(\alpha)$. For such a choice, $\tau(b)$ is nearly second-order \yc{asymptotic }optimal in the sense that for any $\theta \in \Theta_\sigma -\{\theta_0\}$, as $\alpha \rightarrow 0$,
\begin{equation}
\mathbb{E}_{\theta,0}[\tau(b)] - \inf_{T\in C(\alpha)}\mathbb{E}_{\theta,0}[T] =O (\log (\log (1/\alpha))).
\label{eq:optimality1}
\end{equation}
\label{cor:optimalityHP}
\end{Corollary}
The result means that, compared with any procedure (including the optimal procedure) calibrated to have a \yc{false-alarm rate} less than $\alpha$, our procedure incurs an at most $\log(\log(1/\alpha))$ increase in the expected \yc{detection delay}, which is usually a small number. \yc{For instance, even for a conservative case when we set $\alpha=10^{-5}$ to control the false-alarm rate, the number is $\log(\log(1/\alpha))=2.44$.}

\subsection{Sequential change-point detection}

\yc{Now we proceed the proof by leveraging the close connection \cite{lorden1971procedures} between the sequential change-point detection and the one-sided hypothesis test.}
For sequential change-point detection, the two commonly used performance metrics \citep{tartakovsky2014sequential} are the average run length (ARL), denoted by $\mathbb{E}_{\infty}[T]$; and the maximal conditional average delay to detection (CADD), denoted by $\sup_{\nu\geq 0} \mathbb{E}_{\theta, \nu}[T-\nu \mid T> \nu]$. ARL is the expected number of samples between two successive false alarms, and CADD is the expected number of samples needed to detect the change after it occurs. A good procedure should have a large ARL and a small CADD. Similar to the one-sided hypothesis test, one usually choose the threshold large enough so that ARL is larger than a pre-specified level. 

\yc{Similar to Theorem \ref{maintheorem}, we provide an upper bound for the CADD of our ASR and ACM procedures.} 

\begin{Theorem}
Consider the change-point detection procedure $T_{\rm ACM}(b_1)$  in (\ref{ACMprocedure}) and $T_{\rm ASR}(b_2)$  in (\ref{ASRprocedure}). \ycc{For any fixed $k$, let $\{\hat{\theta}_{k,t} \}_{t\geq 1}$ be a sequence of non-anticipating estimators generated by an OCO algorithm $\textsf{a}$.} \yc{Let $b_1=b_2=b$, as $b \rightarrow \infty$} we have that
\begin{equation}
\begin{split}
&\sup_{\nu\geq 0} \mathbb{E}_{\theta, \nu}[T_{\rm ASR}(b)-\nu \mid T_{\rm ASR}(b)> \nu]
\leq 
\sup_{\nu\geq 0} \mathbb{E}_{\theta, \nu}[T_{\rm ACM}(b)-\nu \mid T_{\rm ACM}(b)> \nu] \\ 
\leq &~~
(I(\theta, \theta_0))^{-1}\left(b + \mathbb{E}_{\theta,0}\left[\mathcal R^\textsf{a}_{\tau(b)} \right]+ O(1)  \right).
\end{split} 
\label{eq:CDbound}
\end{equation}
\label{CADD1}
\end{Theorem}

\yc{To prove Theorem \ref{CADD1}, we relate the ASR and ACM procedures to the one-sided hypothesis test and use the fact that when the measure $\mathbb{P}_{\infty}$ is known,
$
\sup_{\nu\geq 0} \mathbb{E}_{\theta, \nu}[T-\nu \mid T> \nu]
$
is attained at $\nu=0$ for both the ASR and the ACM procedures.}
Above, we may apply a similar argument as in Corollary \ref{cor:logregret} to remove the dependence on $\tau(b)$ on the right-hand-side of the inequality. 
We establish the following lower bound for the ARL of the detection procedures, which is needed for proving Corollary \ref{optimality_2}:
\begin{Lemma}[ARL]
Consider the change-point detection procedure $T_{\rm ACM}(b_1)$ in (\ref{ACMprocedure}) and $T_{\rm ASR}(b_2)$ in (\ref{ASRprocedure}). For any fixed $k$, let $\{\hat{\theta}_{k,t} \}_{t\geq 1}$ be any sequence of non-anticipating estimators. Let $b_1=b_2=b$, given a prescribed lower bound $\gamma >0$ for the ARL, we have 
\[
\mathbb{E}_{\infty}[T_{\rm ACM}(b)] \geq \mathbb{E}_{\infty}[T_{\rm ASR}(b)] \geq \gamma,
\] 
provided that $b \geq \log \gamma$.
\label{ARL1}
\end{Lemma}

Lemma \ref{ARL1} shows that given a required lower bound $\gamma$ for ARL, we can choose $b=\log \gamma$ to make the ARL be greater than $\gamma$. This is consistent with earlier works \cite{pollak1987average,lorden2005nonanticipating} which show that the smallest threshold $b$ such that $\mathbb{E}_{\infty}[T_{ACM}(b)] \geq \gamma$ is approximate $\log \gamma$. \yc{However, the bound in Lamma \ref{ARL1} is not tight, since in practice we can set $b = \rho \log\gamma$ for some $\rho \in (0,1)$ to ensure that ARL is greater than $\gamma$.}

Combing the upper bound in Theorem \ref{CADD1} with an existing lower bound for the \yc{CADD} of SRRS procedure in \citep{siegmund2008minimax}, we obtain the following \yc{optimality properties}. 
\begin{Corollary}[Nearly second-order \yc{asymptotic }optimality of ACM and ASR]\label{optimality_2}
\yc{Consider the change-point detection procedure $T_{\rm ACM}(b_1)$  in (\ref{ACMprocedure}) and $T_{\rm ASR}(b_2)$  in (\ref{ASRprocedure}). \ycc{For any fixed $k$, let $\{\hat{\theta}_{k,t} \}_{t\geq 1}$ be a sequence of non-anticipating estimators generated by an OCO algorithm $\textsf{a}$.} Assume that $\mathbb{E}_{\theta,0}[\mathcal R^\textsf{a}_{n}] \leq C\log n$ for any positive integer $n$ and some constant $C>0$. Let $b_1=b_2=b$.} Define $S(\gamma) = \{T: \mathbb{E}_{\infty}[T] \geq \gamma \}$. For $b=\log \gamma$, due to Lemma \ref{ARL1}, both $T_{\rm ASR}(b)$ and $T_{\rm ACM}(b)$ belong to $S(\gamma)$. For such $b$, both $T_{\rm ASR}(b)$ and $T_{\rm ACM}(b)$ are nearly second-order \yc{asymptotic }optimal in the sense that for any $\theta \in \Theta-\{\theta_0\}$ 
\begin{equation}
\begin{split}
& \sup_{\nu\geq 1} \mathbb{E}_{\theta, \nu}[T_{\rm ASR}(b)-\nu+1 \mid T_{\rm ASR}(b)\geq\nu]  \\
& ~~- \inf_{T(b)\in S(\gamma)}\sup_{\nu\geq 1} \mathbb{E}_{\theta, \nu}[T(b)-\nu+1 \mid T(b)\geq\nu] = O(\log\log \gamma ).
\end{split}
\label{eq:optimality2}
\end{equation}
\yc{A} similar expression holds for $T_{\rm ACM}(b)$.
\end{Corollary}
\yc{The result means that, compared with any procedure (including the optimal procedure) calibrated to have a fixed ARL larger than $\gamma$, our procedure incurs an at most $\log(\log \gamma)$ increase in the CADD.}
Comparing (\ref{eq:optimality2}) with (\ref{eq:optimality1}), we note that the ARL $\gamma$ plays the same role as $1/\alpha$ because $1/\gamma$ is roughly the false-alarm rate for sequential change-point detection \cite{lorden1971procedures}. 

\subsection{Example: Regret bound for specific cases}
\label{sec:regretboundexample}
In this subsection, we show that the regret bound $\mathcal{R}_t$ can be expressed as a weighted sum of Bregman divergences between two consecutive estimators. This form of $\mathcal{R}_t$ is useful to show the logarithmic regret for OMD. The following result comes as a modification of \cite{azoury2001relative}.

\begin{Theorem}
Assume that $X_1, X_2, \ldots$ are i.i.d. random variables with density function $f_{\theta}(x)$. Let $\eta_i = 1/i$ in Algorithm \ref{alg1}. Assume that $\{\hat{\theta}_i\}_{i\geq 1}, \{\hat{\mu}_i\}_{i\geq 1}$ are obtained using Algorithm \ref{alg1} and $\hat{\theta}_i = \tilde{\theta}_i$ \yc{(defined in step 7 and 8 of Algorithm \ref{alg1})} for any $i \geq 1$. Then for any $\theta_0 \in \Theta$ and $t\geq 1$,
\[
\mathcal{R}_t = \sum_{i=1}^t i \cdot B_{\Phi^*}(\hat{\mu}_i, \hat{\mu}_{i-1}) =\frac{1}{2} \sum_{i=1}^t i \cdot (\hat{\mu}_i-\hat{\mu}_{i-1})^\intercal [\nabla^2 \Phi^*(\tilde{\mu}_i)](\hat{\mu}_i-\hat{\mu}_{i-1}),
\]
where $\tilde{\mu}_i=\lambda \hat{\mu}_i + (1-\lambda)\hat{\mu}_{i-1}$, for some $\lambda \in (0,1)$. 
\label{thm:regret}
\end{Theorem}

Next, we use Theorem \ref{thm:regret} on a concrete example. The multivariate normal distribution, denoted by $\mathcal{N}(\theta, I_d)$, is parametrized by an unknown mean parameter $\theta$ and a known covariance matrix $I_d$ ($I_d$ is a $d\times d$ identity matrix). Following the notations in subsection \ref{sec:ef}, we know that $\phi(x) = x$, $dH(x) = (1/\sqrt{|2\pi I_d|}) \cdot \exp\left(- x^\intercal  x/2\right)$, $\Theta =\Theta_\sigma = \mathbb{R}^d$ for any $\sigma < 2$, 
$\Phi(\theta) = (1/2) \theta^\intercal \theta$, $\mu=\theta$ and $\Phi^*(\mu) = (1/2) \mu^\intercal \mu$ , where $|\cdot |$ denotes the determinant of a matrix, and $H$ is a probability measure under which the sample follows $\mathcal{N}(0, I_d)$). When the covariance matrix is known to be some $\Sigma \neq I_d$, one can ``whiten'' the vectors by multiplying $\Sigma^{-1/2}$ to obtain the situation here.
\begin{Corollary}[Upper bound for the expected regret, Gaussian]
Assume $X_1, X_2, \ldots$ are i.i.d. following $\mathcal{N}(\theta, I_d)$ with some  $\theta \in \mathbb{R}^d$. Assume that $\{\hat{\theta}_i\}_{i\geq 1}, \{\hat{\mu}_i \}_{i\geq 1}$ are obtained using Algorithm \ref{alg1} with $\eta_i = 1/i$ and $\Gamma =\mathbb{R}^d$. For any $t>0$, we have that for some constant $C_1>0$ that depends on $\theta$,
\[
\mathbb{E}_{\theta,0}[\mathcal{R}_t] \leq C_1d\log t/2. 
\]
\label{cor:regret_normal}
\end{Corollary}

The following calculations justify Corollary \ref{cor:regret_normal}, which also serve as an example of how to use regret bound. First, the assumption $\hat{\theta}_t = \tilde{\theta}_t$ in Theorem \ref{thm:regret} is satisfied for the following reasons. Consider $\Gamma = \mathbb{R}^d$ is the full space. According to Algorithm \ref{alg1}, using the non-negativity of the Bregman divergence, we have 
$
\hat{\theta}_t = \mathop{\arg\min}_{u \in \Gamma} B_{\Phi}(u, \tilde{\theta}_t) = \tilde{\theta}_t.
$
Then the regret bound can be written as
\begin{equation}
\begin{split}
\mathcal{R}_t =& \frac{1}{2} (\hat{\mu}_1 - \hat{\mu}_0)^{\intercal}(\hat{\mu}_1 - \hat{\mu}_0) + \frac12 \sum_{i=2}^t [ i \cdot(\hat{\mu}_i-\hat{\mu}_{i-1})^\intercal (\hat{\mu}_i-\hat{\mu}_{i-1}) ] \\
=& \frac{1}{2} (X_1 - \theta_0)^{\intercal}(X_1 - \theta_0) +  \frac12 \sum_{i=2}^t (\hat{\mu}_i-\hat{\mu}_{i-1})^\intercal (\phi(X_i) - \hat{\mu}_{i-1}).
\end{split}\nonumber
\end{equation}
Since the step-size $\eta_i = 1/i$, the second term in the above equation can be written as: 
\begin{equation} 
\begin{split}
&\frac{1}{2} \sum_{i=2}^t (\hat{\mu}_i-\hat{\mu}_{i-1})^\intercal (\phi(X_i) - \hat{\mu}_{i-1}) \\
=&\frac{1}{2} \sum_{i=2}^t (\hat{\mu}_i-\hat{\mu}_{i-1})^\intercal (\phi(X_i) + \hat{\mu}_{i}) - \sum_{i=2}^t \frac{1}{2}(\hat{\mu}_i-\hat{\mu}_{i-1})^\intercal (\hat{\mu}_{i-1} + \hat{\mu}_{i}) \\
=& \sum_{i=2}^t \frac{1}{2(i-1)}(\phi(X_i) - \hat{\mu}_{i})^\intercal (\phi(X_i) + \hat{\mu}_{i}) + \sum_{i=2}^t \frac{1}{2}(\left\Vert\hat{\mu}_{i-1}\right\Vert^2 - \left\Vert \hat{\mu}_i\right\Vert^2) \\
=& \sum_{i=2}^t \frac{1}{2(i-1)} \left\Vert X_i\right\Vert^2 - \sum_{i=2}^t \frac{1}{2(i-1)} \left\Vert \hat{\mu}_i \right\Vert^2  + \frac{1}{2}\left\Vert\hat{\mu}_1\right\Vert^2 - \frac{1}{2}\left\Vert\hat{\mu}_t\right\Vert^2. 
\end{split} \nonumber
\end{equation}
Combining above, we have
\[
\mathbb{E}_{\theta,0}[\mathcal{R}_t] \leq \frac{1}{2} \mathbb{E}_{\theta,0}[(X_1 - \theta_0)^{\intercal}(X_1 - \theta_0)] + \frac{1}{2} \sum_{i=2}^t \frac{1}{i-1} \mathbb{E}_{\theta,0}[\left\Vert X_i\right\Vert^2] +\frac{1}{2} \mathbb{E}_{\theta,0}[\left\Vert X_1\right\Vert^2].
\]
Finally, since $\mathbb{E}_{\theta,0}[\left\Vert X_i\right\Vert^2] = d(1+\theta^2)$ for any $i\geq 1$, we obtain desired result. Thus, with i.i.d. multivariate normal samples, the expected regret grows logarithmically with the number of samples. 

Using the similar calculations, we can also bound the expected regret in the general case. As shown in the proof above for Corollary \ref{cor:regret_normal}, the dominating term for $\mathcal{R}_t$ can be rewritten as 
\[
\sum_{i=2}^t \frac{1}{2(i-1)}(\phi(X_i) - \hat{\mu}_{i})^\intercal [\nabla^2 \Phi^*(\tilde{\mu}_i)](\phi(X_i) + \hat{\mu}_{i}),
\]
where $\tilde{\mu}_i$ is a convex combination of $\hat{\mu}_{i-1}$ and $\hat{\mu}_i$. For an arbitrary distribution, the term $(\phi(X_i) - \hat{\mu}_{i})^\intercal [\nabla^2 \Phi^*(\tilde{\mu}_i)](\phi(X_i) + \hat{\mu}_{i})$ can be viewed as a local normal distribution with the changing curvature $\nabla^2 \Phi^*(\tilde{\mu}_i)$. Thus, it is possible to prove case-by-case the $O(\log t)$-style bounds \ycc{by making more assumptions about the distributions}. \ycc{Recall the notation $\Theta_\sigma$ in subsection \ref{sec:ef} such that $-\log f_{\theta}(x)$ is $\sigma$-strongly convex over $\Theta_\sigma$. Let $\|\cdot\|_2$ denote the $\ell_2$ norm. Moreover, we assume that the true parameter belongs to a set $\Gamma$ that is a closed and convex subset of $\Theta_\sigma$ such that $\sup_{\theta \in \Gamma} \| \nabla \Phi(\theta)\|_2  \leq M$ for some constant $M$. Thus, one can show that $-\log f_{\theta}(x)$ is not only $\sigma$-strongly convex but also $M$-strongly smooth over $\Gamma$.
Theorem 3 in \cite{raginsky2012sequential} shows that for all $\theta \in \Gamma$ and $n\geq 1$, consider that $\{\hat{\theta}_i\}_{i\geq 1}$ is obtained by OMD, then 
\[
\mathbb{E}_{\theta,0}[\mathcal{R}_n] \leq \frac{\mathbb{E}_{\theta,0}\left[\left(\frac{1}{2}\max_{1\leq i\leq n} \|X_i\|_2 + \frac{1}{2}M \right)^2 \right]}{\sigma} \cdot (\log n + 1).
\]
Therefore, for any bounded distributions within the exponential family, we achieve a logarithmic regret. This logarithmic regret is valid for Bernoulli distribution, Beta distribution and some truncated versions of classic distributions (e.g., truncated Gaussian distribution, truncated Gamma distribution and truncated Geometric distribution analyzed in \cite{alqanoo2014on}). 
}

\section{Numerical examples}
\label{sec:simulation}

In this section, we present some synthetic examples to demonstrate the good performance of our methods. We will focus on ACM and ASR for sequential change-point detection. \yc{In the following, we consider the window-limited versions (see Remark \ref{remark_window}) of ACM and ASR with window size $w=100$.} \yc{Recall that when the measure $\mathbb{P}_{\infty}$ is known,
$
\sup_{\nu\geq 0} \mathbb{E}_{\theta, \nu}[T-\nu \mid T> \nu]
$
is attained at $\nu=0$ for both ASR and ACM procedures (a proof can be found in the proof of Theorem \ref{CADD1}). Therefore, in the following experiments we define the expected detection delay (EDD) as $\mathbb{E}_{\theta, 0}[T]$  for a stopping time $T$. To compare the performance between different detection procedures, we determine the threshold for each detection procedure by Monte-Carlo simulations such that the ARL for each procedure is about $10000$. Below, we denote $\left\Vert\cdot\right\Vert_2$, $\left\Vert\cdot\right\Vert_1$ and $\left\Vert\cdot\right\Vert_0$ as the $\ell_2$ norm, $\ell_1$ norm and $\ell_0$ norm defined as the number of non-zero entries, respectively. The following experiments are all run on the same Macbook Air with an Intel i7 Core CPU.}

\subsection{Detecting sparse mean-shift of multivariate normal distribution}

We consider detect the \yc{sparse mean shift for} multivariate normal distribution. \yc{Specifically, we assume that the pre-change distribution is $\mathcal{N}(0,I_d)$ and the post-change distribution is $\mathcal{N}(\theta, I_d)$ for some unknown $\theta \in \{ \theta\in \mathbb{R}^d: \|\theta\|_0 \leq s\}$, where $s$ is called the \textit{sparsity} of the mean shift.} Sparse mean shift detection is of particular interest in sensor networks \cite{xie2013sequential, siegmund2011detecting}. 
For this Gaussian case, the Bregman divergence is given by 
$
B_{\Phi}(\theta_1, \theta_2) = I(\theta_2, \theta_1) =  \| \theta_1 - \theta_2 \|_2^2/2.
$ 
Therefore, the projection onto $\Gamma$ in Algorithm \ref{alg1} is a Euclidean projection onto a convex set, which in many cases can be implemented efficiently. As a frequently used convex relaxation of the $\ell_0$-norm ball, we set $\Gamma = \{\theta: \| \theta \|_1 \leq s \}$ (\yao{it is known that imposing an $\ell_1$ constraint leads to sparse solution; see, e.g.,  [48]}). Then, the projection onto $\ell_1$ ball can be computed very efficiently via a simple soft-thresholding technique \cite{duchi2008efficient}. 

Two benchmark procedures are the CUSUM and the GLR. For \yc{the} CUSUM procedure, we specify a nominal post-change mean, which is an all-one vector. \ycc{If  knowing the post-change mean is sparse, we can also use the shrinkage estimator presented in \cite{wang2015large}, which performs hard or soft thresholding of the estimated post-change mean parameter.
}
Our procedures are $T_{ASR}(b)$ and $T_{ACM}(b)$ with $\Gamma = \mathbb{R}^d$ and $\Gamma = \{\theta: \| \theta \|_1 \leq 5 \}$.  In the following experiments, we run $10000$ Monte Carlo trials to obtain each simulated EDD. 

In the experiments, we set $d=20$. The post-change distributions are $\mathcal{N}(\theta, I_d)$, where $100p\%$ entry of $\theta$ is $1$ and others are $0$, and the location of nonzero entries are random. Table \ref{tab:normal_mean_p} shows the EDDs versus the proportion $p$. Note that our procedures incur little performance loss compared with the GLR procedure and the CUSUM procedure. Notably, $T_{ACM}(b)$ with $\Gamma =  \{\theta: \| \theta \|_1 \leq 5 \}$ performs almost the same as the GLR procedure and much better than the CUSUM procedure when $p$ is small. This shows the advantage of projection when the true parameter is sparse.


\begin{table}[h]
\centering
\begin{tabular}{|c|c|c|c|c|c|c|}
\hline
 & $p=0.1$  & $p= 0.2$  &  $p=0.3$  & $p=0.4$ &  $p=0.5$ & $p=0.6$ \\ 
\hline
CUSUM & 188.60 & 146.45 &  64.30 &  18.97  &  7.18  &  3.77 \\
\hline
Shrinkage & 17.19 &    9.25 & 6.38 & 4.96 &  4.07 &  3.55 \\
\hline
GLR & 19.10 &  10.09 &   7.00  &  5.49 &   4.50   & 3.86 \\ 
\hline
ASR &  45.22  & 19.55 &  12.62  &  8.90  &  7.02 &   5.90 \\ 
\hline
ACM & 45.60 &  19.93  & 12.50  &  9.00   & 7.03  &  5.87 \\
\hline
ASR-$\ell$1 & 45.81  & 19.94 &  12.45  &  8.92 &   6.97  &  5.89 \\
\hline
ACM-$\ell$1 & 19.24  & 10.17  &  7.51  &   6.11  &  5.41  &  4.92 \\
\hline
\end{tabular}
\caption{Comparison of \yc{the EDDs in detecting the sparse mean shift of multivariate Gaussian distribution.} Below, ``CUSUM'': CUSUM procedure with pre-specified all-one vector as post-change parameter; ``Shrinkage'': component-wise shrinkage estimator in \cite{wang2015large}; ``GLR'': GLR procedure; ``ASR'': $T_{\rm ASR}(b)$ with $\Gamma=\mathbb{R}^d$; ``ACM'': $T_{\rm ACM}(b)$ with $\Gamma=\mathbb{R}^d$; ``ASR-L1'': $T_{\rm ASR}(b)$ with $\Gamma= \{\theta: \| \theta \|_1 \leq 5 \}$; ``ACM-L1'': $T_{ACM}(b)$ with $\Gamma= \{\theta: \| \theta \|_1 \leq 5 \}$. $p$ is the proportion of non-zero entries in $\theta$. We run $10000$ Monte Carlo trials to obtain each value. For each \yc{value}, the standard deviation is less than one half of the \yc{value}.}
\label{tab:normal_mean_p}
\end{table}

\subsection{Detecting the scale change in Gamma distribution}
We consider an example that detects the scale change in Gamma distributions. Assume that we observe a sequence $X_1, X_2 \ldots$ of samples drawn from $\mbox{Gamma}(\alpha, \beta)$ for some $\alpha, \beta>0$, with the probability density function given by $f_{\alpha, \beta}(x) = \exp(-x\beta)x^{\alpha-1}\beta^{\alpha}/\tilde{\Gamma}(\alpha)$ (to avoid confusion with the $\Gamma$ parameter in Algorithm \ref{alg1} we use $\tilde{\Gamma}(\cdot)$ to denote the Gamma function). The parameter $\alpha^{-1}$ is called the dispersion parameter that scales the loss and the divergences. For simplicity, we fix $\alpha =1$ just like we fix the variance in the Gaussian case. The specifications in the Algorthm \ref{alg1} are as follows: $\theta = -\beta$, $\Theta = (-\infty, 0)$, $\phi(x) = x$, $dH(x) = 1$, $\Phi(\theta) = -\log (-\theta)$, $\mu = -1/\theta$ and $\Phi^*(\mu) = -1-\log \mu$.
Assume that the pre-change distribution is $\mbox{Gamma}(1,1)$ and the post-change distribution is $\mbox{Gamma}(1,\beta)$ for some unknown $\beta>0$. 
\ycc{We compare our algorithms with CUSUM, GLR and non-ancitipating estimator based on the method of moment (MOM) estimator in \cite{lorden2005nonanticipating}. For the CUSUM procedure, we specify the post-change $\beta$ to be $2$. 
} The results are shown in Table \ref{tab:gamma}.  CUSUM fails to detect the change when $\beta=0.1$, which is far away from the pre-specified post-change parameter $\beta=2$. We can see that performance loss of the proposed ACM method compared with GLR and MOM is very small.

\begin{table}[h]
\centering
\begin{tabular}{|c|c|c|c|c|c|}
\hline
 & $\beta = 0.1$  & $\beta = 0.5$  &  $\beta = 2$  & $\beta = 5$ &  $\beta = 10$\\ 
\hline
CUSUM & NaN & 481.2 & 33.75 &  14.37  &  12.04 \\
\hline
MOM & 3.41 &  32.87 & 40.86 &   11.42 &  7.21\\
\hline
GLR & 2.40 &  23.79 &   33.29  &  9.07 &   5.67  \\ 
\hline
ASR &  3.95  & 32.34 &  45.18  &  13.45 &  8.55 \\ 
\hline
ACM & 3.70 &  31.80  & 47.20  &  12.42  & 7.87 \\
\hline
\end{tabular}
\caption{\yc{Comparison of the EDDs in detecting the scale change in Gamma distribution. Below, ``CUSUM'': CUSUM procedure with pre-specified post-change parameter $\beta=2$; ``MOM'': Method of Moments estimator method; ``GLR'': GLR procedure; ``ASR'': $T_{\rm ASR}(b)$ with $\Gamma=(-\infty,0)$; ``ACM'': $T_{\rm ACM}(b)$ with $\Gamma=(-\infty,0)$. We run $10000$ Monte Carlo trials to obtain each value. For each value, the standard deviation is less than one half of the value.}}
\label{tab:gamma}
\end{table}

\subsection{Communication-rate change detection with Erd\H{o}s-R\'enyi model }

Next, we consider a problem to detect the communication-rate change in a network, which is a model for social network data. Suppose we observe communication between nodes in a network over time, represented as a sequence of (symmetric) adjacency matrices of the network. At time $t$, if node $i$ and node $j$ communicates, then the adjacency matrix has 1 on the $ij$th and $ji$th entries (thus it forms an undirected graph). The nodes that do not communicate have 0 on the corresponding entries. We model such communication patterns using the Erdos-Renyi random graph model. Each edge has a fixed probability of being present or absent, independently of the other edges. Under the null hypothesis, each edge is a Bernoulli random variable that takes values $1$ with known probability $p$ and value $0$ with probability $1-p$. Under the alternative hypothesis, there exists an unknown time $\kappa$, after which a small subset of edges occur with an unknown and different probability $p' \neq p$.

In the experiments, we set $N=20$ and $d=190$. For the pre-change parameters, we set $p_i=0.2$ for all $i=1,\ldots,d$. For the post-change parameters, we randomly select $n$ out of the $190$ edges, denoted by $\mathcal{E}$, and set $p_i = 0.8$ for $i \in \mathcal{E}$ and $p_i = 0.2$ for $i \notin \mathcal{E}$. As said before, let the change happen at time $\nu=0$ (since the upper bound for EDD is achieved at $\nu=0$ as argued in the proof of Theorem \ref{CADD1}). To implement CUSUM, we specify the post-change parameters $p_i = 0.8$ for all $i=1,\ldots,d$. 

The results are shown in Table \ref{tab:ERgraph}. Our procedures are better than CUSUM procedure when $n$ is small since the post-change parameters used in CUSUM procedure is far from the true parameter. Compared with GLR procedure, our methods have a small performance loss, and the loss is almost negligible as $n$ approaches to $d=190$. 


\begin{table}[h]
\centering
\begin{tabular}{|c|c|c|c|c|c|c|c|}
\hline
 & $n= 78$  &  $n=100$  & $n=120$ &  $n=150$ & $n=170$ & $n=190$\\
\hline
CUSUM & 473.11   & 2.06 &   2.00   & 2.00  &  2.00   & 2.00\\
\hline
GLR &   2.00  &  1.96  &  1.27 &   1.00  &  1.00  &  1.00 \\
\hline
ASR  &  8.64 &   6.39  &  5.08 &   3.92  &  3.36  &  2.94\\ 
\hline
ACM &  8.67 &   6.37  &  5.07   & 3.88  &  3.32  &  2.94 \\
\hline
\end{tabular}
\caption{Comparison of the EDDs in detecting the changes of the communication-rates in a network.  Below, ``CUSUM'': CUSUM procedure with pre-specified post-change parameters $p=0.8$ ; ``GLR'': GLR procedure; ``ASR'': $T_{\rm ASR}(b)$ with $\Gamma=\mathbb{R}$; ``ACM'': $T_{\rm ACM}(b)$ with $\Gamma=\mathbb{R}$. We run $10000$ Monte Carlo trials to obtain each value. For each value, the standard deviation is less than one half of the value.}
\label{tab:ERgraph}
\end{table}

Below are the specifications of Algorithm \ref{alg1} in this case. For Bernoulli distribution with unknown parameter $p$, the natural parameter $\theta$ is equal to $\log(p/(1-p))$. Thus, we have $\Theta = \mathbb{R}$, $\phi(x) = x$, $dH(x) = 1$, $\Phi(\theta) = \log(1+\exp(\theta))$, $\mu=\exp(\theta)/(1+\exp(\theta))$ and $\Phi^*(\mu) = \mu \log \mu + (1-\mu)\log(1-\mu)$.

\subsection{\yao{Point process change-point detection: Poisson to Hawkes processes}}\label{sec:poisson_hawkes}

\yao{
In this example, to illustrate the situation in Section \ref{sec:social_net}, we consider a case where a homogeneous Poisson process switches to a Hawkes process (see, e.g., \cite{LiXie17}); this can be viewed as a simplest case in Section \ref{sec:social_net} with one node. We construct ACM and ASR procedures. In this case, the MLE for the unknown post-change parameter cannot be found in close-form, yet ACM and ASR can be easily constructed and give reasonably good performance, although our theory no longer holds in this case due to the lack of i.i.d. samples.   }

\yao{The Hawkes process can be viewed as a non-homogeneous Poisson process where the intensity is influenced by historical events. The data consist of a sequence of events occurring at times $\{t_1, t_2, \ldots, t_n\}$ before a time horizon $T$: $t_i\leq T$. Assume the intensity of the Poisson process is $\lambda_s, s\in (0, T)$ and there may exists a change-point $\kappa \in (0, T)$ such that the process changes. The null and alternative hypothesis tests are
\begin{equation*} \label{test_poi_haw_1d}
\left\{
\begin{array}{ll}
 \textsf{H}_0: &  \lambda_s =\mu, \quad 0 < s <T;   \\
 \textsf{H}_1: &  \lambda _s=\mu, \quad 0 < s <\kappa,\\
&   \lambda_s=\mu+\theta  \sum_{\kappa<t_j<s} \varphi(s-t_j), \quad \kappa < s<T, 
\end{array}
\right.
\end{equation*}
where $\mu$ is a known baseline intensity, $\theta >0$ is unknown magnitude of the change, $\varphi(s)=\beta e^{-\beta s}$ is the normalized kernel function with pre-specified parameter $\beta>0$, which captures the influence from the past events. We treat the post-change influence parameter $\theta$ as unknown  since it represents an anomaly. }

\yao{We first use a sliding window to convert the event times into a sequence of vectors with overlapping events. Assume of size of the sliding window is $L$. For a given scanning time $T_i \leq T$, we map all the events in $[T_i-L, T_i]$ to a vector $X_i = [t_{(1)}, \ldots, t_{(m_i)}]^\intercal$, $t_{(i)} \in [T_i-L, T_i]$, where $m_i$ is the number of events falling into the window. Note that $X_i$ can have different length for different $i$. Consider a set of scanning times $T_1, T_2, \ldots, T_t$. This maps the event times into a sequence of vectors $X_1, X_2, \ldots, X_t$ of lengthes $m_1$, $m_2$, $\ldots$, $m_t$. These scanning times can be arbitrary; here we set them to be event times so that there are at least one sample per sliding window. }

\yao{For a hypothetical change-point location $k$, it can be shown that the log-likelihood ratio (between the Hawkes process and the Poisson process) as a function of $\theta$, is given by
\begin{equation}
\ell(\theta|X_i)=\sum_{t_q \in (T_i-L, T_i)}  \mbox{log} \left[ \mu +  \theta \sum_{t_j \in (T_i-L, t_q)} \beta e^{-\beta (t_q-t_j)} \right]  - \mu L  - \theta \sum_{t_q \in (T_i-L, T_i)} \left[1 -e^{-\beta (T_i-t_q)}  \right].
\label{loglikelihood_poi}
\end{equation}
Now based on this sliding window approach, we can approximate the original change-point detection problem as the following. Without change, $X_1, \ldots, X_t$ are sampled from a Poisson process. 
Under the alternative, the change occurs at some time such that $X_1, \ldots, X_\kappa$ are sampled from a Poisson process, and $X_{\kappa+1}, \ldots, X_t$ are sampled from a Hawkes process with parameter $\theta$, rather than a Poisson process. We
 define the estimator of $\theta$, for assumed change-point location $\kappa = k$ as follows
\begin{equation}
\hat \theta_{k, i} \triangleq \hat \theta_{k,i}(X_k, \ldots, X_i) 
= \hat \theta_{k,i}(t_\ell\in [T_k, T_i])
\end{equation}
Now, consider $k \in [i-w, i-1]$, and keep $w$ estimators: $\hat \theta_{i-w,i}, \ldots, \hat \theta_{i-1,i}$. The update for each estimator is based on stochastic gradient descent. By taking derivative with respect to $\theta$, we have
\[
\frac{\partial \ell(\theta|X_i)}{\partial \alpha}  = \sum\limits_{t_q \in (T_i-L,T_i)} \frac{ \sum_{t_j \in (T_i-L, t_q)} \beta e^{-\beta (t_q-t_j)} }{ \mu +  \theta \sum_{t_j \in (T_i-L, t_q)} \beta e^{-\beta (t_q-t_j)}}  - \sum_{t_q \in (T_i-L, T_i)} \left[1 -e^{-\beta (T_i-t_q)} \right],
\]
Note that there is no close form expression for the MLE, which the solution to the above equation. We perform stochastic gradient descent instead
\[
\hat\theta_{k,i+1} = \hat\theta_{k,i} - \gamma \frac{\partial \ell(\theta|X_{i+1}) }{\partial \theta}\Big|_{\theta = \hat\theta_{k,i}}, \quad k = i-w+1, i-w, \ldots, i,
\]
where $\gamma >0$ is the step-size.
Now we can apply the ACM and ASR procedures, by using the fact that $f_{\hat \theta_{k,t}}(X_{t+1})/f_{\theta_0}(X_{t+1})  = \ell(\hat\theta_{k,t}|X_{t+1})$ and calculating using (\ref{loglikelihood_poi}).}

\yao{Table. \ref{tab:EDD_A_hawkes} shows the EDD for different $\alpha$. Here we choose the threshold such that ARL is 5000. We see that the scheme has a reasonably good performance, the detection delay decreases as the true signal strength $\theta$ increases. }

\begin{table}[h]
\centering
\begin{tabular}{|c|c|c|c|c|c|}
\hline
   & $\theta = 0.4$ &  $\theta = 0.5$ & $\theta = 0.5$ & $\theta = 0.7$\\
\hline
 ACM    & 33.03 &  27.75 &   20.39  & 16.16 \\
 \hline
 ASR   & 38.59  & 24.96 &  20.17 &  13.91 \\
\hline
\end{tabular}
\caption{Point process change-point detection: EDD of ACM and ASR procedures for various values of true $\theta$; ARL of the procedure is controlled to be 5000 by selecting threshold via Monte Carlo simulation.}
\label{tab:EDD_A_hawkes}
\end{table}

%

\section{Conclusion}
\label{sec:conclusion}

In this paper, we consider sequential hypothesis testing and change-point detection with computationally efficient one-sample update schemes obtained from online mirror descent. We show that the loss of the statistical efficiency caused by the online mirror descent estimator (replacing the exact maximum likelihood estimator using the complete historical data) is related to the regret incurred by the online convex optimization procedure. The result can be generalized to any estimation method with logarithmic regret bound. This result sheds lights on the relationship between the statistical detection procedures and the online convex optimization.

\acknowledgments{This research was supported in part by National Science Foundation (NSF) NSF CCF-1442635, CMMI-1538746, NSF CAREER CCF-1650913 to Yao Xie. \yao{We would like to thank the anonymous reviewers to provide insightful comments.}}

\authorcontributions{Yang Cao, Yao Xie, and Huan Xu conceived the idea and performed the theoretical part of the paper; Liyan Xie helped with numerical examples of the manuscript.}

%

\bibliography{ref_SCD_Tradeoff}

\appendix

\section{Proofs}

\begin{proof}[Proof of Theorem \ref{maintheorem}]
In the proof, for the simplicity of notation we use $N$ to denote $\tau(b)$. Recall $\theta$ is the true parameter.
Define that 
\[
S_t^{\theta} = \sum_{i=1}^t \log \frac{f_{\theta}(X_i)}{f_{\theta_0}(X_i)}. 
\] Then under the measure $\mathbb{P}_{\theta,0}$, $S_t$ is a random walk with i.i.d. increment. Then, by Wald's identity (e.g., \cite{siegmund1985sequential}) we have that 
\begin{equation}
\mathbb{E}_{\theta,0}[S_N^{\theta}] = \mathbb{E}_{\theta,0}[N] \cdot I(\theta,\theta_0).
\label{waldequation}
\end{equation}

On the other hand, let $\theta_N^*$ denote the MLE based on $(X_1, \ldots, X_N)$. The key to the proof is to decompose the stopped process $S_N^{\theta}$ as a summation of three terms as follows:
\begin{equation}
S_N^{\theta} = \sum_{i=1}^N \log\frac{f_{\theta}(X_i)}{f_{\theta_N^*}(X_i)} + \sum_{i=1}^N \log \frac{f_{\theta_N^*}(X_i)}{f_{\hat{\theta}_{i-1}}(X_i)} + \sum_{i=1}^N \log \frac{f_{\hat{\theta}_{i-1}}(X_i)}{f_{\theta_0}(X_i)},
\label{decomposition}
\end{equation}
Note that the first term of the decomposition on the right-hand side of (\ref{decomposition}) is always non-positive since
\[
\sum_{i=1}^N \log\frac{f_{\theta}(X_i)}{f_{\theta_N^*}(X_i)}=\sum_{i=1}^N \log f_{\theta}(X_i) - \sup_{\tilde{\theta} \in \Theta} \sum_{i=1}^N \log f_{\tilde{\theta}}(X_i) \leq 0.
\]
Therefore we have
\[
\mathbb{E}_{\theta,0}[S_N^{\theta}] \leq \mathbb{E}_{\theta,0}\left[\sum_{i=1}^N \log \frac{f_{\theta_N^*}(X_i)}{f_{\hat{\theta}_{i-1}}(X_i)}\right] + \mathbb{E}_{\theta,0}\left[\sum_{i=1}^N \log \frac{f_{\hat{\theta}_{i-1}}(X_i)}{f_{\theta_0}(X_i)}\right].
\]

Now consider the third term in the decomposition (\ref{decomposition}). Similar to the proof of equation (5.109) in \cite{tartakovsky2014sequential}, we claim that its expectation under measure $\mathbb{P}_{\theta, 0}$ is upper bounded by $b/I(\theta, \theta_0)+O(1)$ as $b\rightarrow \infty$. Next, we prove the claim. For any positive integer $n$, we further decompose the third term in (\ref{decomposition}) as 
\begin{equation}
\sum_{i=1}^n \log \frac{f_{\hat{\theta}_{i-1}}(X_i)}{f_{\theta_0}(X_i)} = M_n(\theta) - G_n(\theta) + m_n(\theta, \theta_0) + nI(\theta, \theta_0),
\label{decomposition2}
\end{equation}
where 
\[
M_n(\theta) = \sum_{i=1}^n \log \frac{ f_{\hat{\theta}_{i-1}}(X_i)}{ f_{\theta}(X_i)} + G_n(\theta),
\]
\[
G_n(\theta) = \sum_{i=1}^n I(\theta, \hat{\theta}_{i-1}),
\]
and
\[
m_n(\theta, \theta_0) = \sum_{i=1}^n \log \frac{ f_{\theta}(X_i)}{f_{\theta_0}(X_i)} - nI(\theta, \theta_0).
\]
The decomposition of (\ref{decomposition2}) consists of stochastic processes $\{M_n(\theta)\}$ and $\{m_n(\theta, \theta_0)\}$, which are both $\mathbb{P}_{\theta,0}$-martingales with zero expectation, i.e., $\mathbb{E}_{\theta,0}[M_n(\theta)] = \mathbb{E}_{\theta,0}[m_n(\theta, \theta_0)] = 0$ for any positive integer $n$. Since for exponential family, the log-partition function $\Phi(\theta)$ is bounded, by the inequalities for martingales \cite{lipster1989theory} we have that 
\begin{equation}
\mathbb{E}_{\theta,0}|M_n(\theta)| \leq C_1 \sqrt{n}, \quad 
\mathbb{E}_{\theta,0}|m_n(\theta, \theta_0)| \leq C_2 \sqrt{n},
\label{martingaleinequality}
\end{equation}
where $C_1$ and $C_2$ are two absolute constants that do not depend on $n$. \ycc{Moreover, we observe that for all $\theta \in \Theta$, 
\[
\mathbb{E}_{\theta,0}[G_n(\theta)] \leq \mathbb{E}_{\theta,0}\left[\max_{\tilde{\theta} \in \Theta} G_n(\tilde{\theta})\right] = \mathbb{E}_{\theta,0}[\mathcal{R}_n(\theta)] \leq C\log n.
\]}
Therefore, applying (\ref{martingaleinequality}), we have that $n^{-1}G_n(\theta), n^{-1}M_n(\theta)$ and $n^{-1}m_n(\theta, \theta_0)$ converge to $0$ almost surely. Moreover, the convergence is $\mathbb{P}_{\theta,0}$-$r$-quickly for $r=1$. \yc{We say that $n^{-1} A_n$ converges $\mathbb{P}_{\theta,0}$-$r$-quickly to a constant $I$ if $\mathbb{E}_{\theta,0}[\mathcal G(\epsilon)]^r < \infty$ for all $\epsilon>0$, where $\mathcal G(\epsilon) = \mbox{sup}\{n\geq 1: |n^{-1}A_n - I| > \epsilon\}$ is the last time when $n^{-1}A_n$ leaves the interval $[I-\epsilon, I+\epsilon]$ (for more details, we refer the readers to Section 2.4.3 of \cite{tartakovsky2014sequential}).}
Therefore, dividing both sides of (\ref{decomposition2}) by n, we obtain $n^{-1} \sum_{i=1}^n \log (f_{\hat{\theta}_{i-1}}(X_i)/ f_{\theta_0}(X_i))$ converges \yc{$\mathbb{P}_{\theta,0}$-}$1$-quickly to $I(\theta, \theta_0)$. 

For $\epsilon>0$, we now define the last entry time
\[
L(\epsilon) = \sup\left\{ n\geq 1: \left| \frac{1}{I(\theta, \theta_0)}\sum_{i=1}^n \log \frac{f_{\hat{\theta}_{i-1}}(X_i)}{f_{\theta_0}(X_i)} - n \right| > \epsilon n \right\}.
\]
By the definition of \yc{$\mathbb{P}_{\theta,0}$-}$1$-quickly convergence \yc{and the finiteness of $I(\theta, \theta_0)$}, we have that $\mathbb{E}_{\theta,0}[L(\epsilon)]<+\infty$ for all $\epsilon>0$. In the following, define a scaled threshold $\tilde{b} = b/I(\theta, \theta_0)$. Observe that conditioning on the event $\{L(\epsilon)+1 < N < +\infty \}$, we have that
\[
(1-\epsilon)(N-1)I(\theta, \theta_0) < \sum_{i=1}^{N-1} \log \frac{f_{\hat{\theta}_{i-1}}(X_i)}{f_{\theta_0}(X_i)} < b.
\]
Therefore, conditioning on the event $\{L(\epsilon)+1 < N < +\infty \}$,
we have that
$ N < 1 + b/(1-\epsilon)$.
Hence, for any $0<\epsilon<1$, we have
\begin{equation}
N \leq 1 + \mathbb{I}(\{N>L(\epsilon)+1\}) \cdot \frac{\tilde{b}}{1-\epsilon} + \mathbb{I}(\{N\leq L(\epsilon)+1 \}) \cdot L(\epsilon) \leq 1+\frac{\tilde{b}}{1-\epsilon}+L(\epsilon).
\label{eq:conditioninequality}
\end{equation}
Since $\mathbb{E}_{\theta,0}[L(\epsilon)]<\infty$ for any $\epsilon>0$, from (\ref{eq:conditioninequality}) above, we have that the third term in (\ref{decomposition}) is upper bounded by $\tilde{b}+\textsf{O}(1)$.

Finally, the second term in (\ref{decomposition}) can be written as
\[
\sum_{i=1}^N \log \frac{f_{\theta_N^*}(X_i)}{f_{\hat{\theta}_{i-1}}(X_i)} = 
\sum_{i=1}^N -\log f_{\hat{\theta}_{i-1}}(X_i) - \inf_{\tilde{\theta} \in \Theta} \sum_{i=1}^N -\log f_{\tilde{\theta}}(X_i),
\]
which is just the regret defined in (\ref{regret}) for the online estimators: $\mathcal R_t$, when the loss function is defined to be the negative likelihood function. Then, the theorem is proven by combining the above analysis for the three terms in (\ref{decomposition}) and (\ref{waldequation}).
\end{proof}

\begin{proof}[Proof of Corollary \ref{cor:logregret}]

\ycc{First, we can relate the expected regret at the stopping time to the expected stopping time, using the following chain of equalities and inequalities 
\begin{equation}
\mathbb{E}_{\theta,0}[\mathcal{R}_{\tau(b)}] = \mathbb{E}_{\theta,0}[\mathbb{E}_{\theta,0}[\mathcal{R}_n \mid \tau(b)=n]]
\leq  \mathbb{E}_{\theta,0}[C\log \tau(b)] 
\leq  C \log \mathbb{E}_{\theta,0} [\tau(b)],
\label{jensen1}
\end{equation}}
where the first equality uses iterative expectation, the first inequality uses the assumption of the logarithmic regret in the statement of Corollary \ref{cor:logregret}, and the second inequality is due to Jensen's inequality. 
Let $\alpha = (b+O(1))/I(\theta, \theta_0)$, $\beta = C/I(\theta, \theta_0)$ and $x = \mathbb{E}_{\theta,0}[\tau(b)]$. Applying (\ref{jensen1}), the upper bound in equation (\ref{sequentialHPdelay}) becomes $x\leq \alpha + \beta \log(x)$. From this, we have $x \leq O(\alpha)$. Taking logarithm on both sides and using the fact that $\max\{a_1+a_2\} \leq a_1+a_2 \leq 2\max\{a_1, a_2\}$ for $a_1, a_2\geq0$, $\log(x) \leq \max\{ \log (2\alpha), \log (2\beta \log x) \} \leq \log(\alpha) + o(\log b)$. Therefore, we have that $x\leq \alpha+\beta (\log(\alpha) + o(\log b))$. Using this argument, we obtain 
\begin{equation}
\mathbb{E}_{\theta,0}[\tau(b)] \leq \frac{b}{I(\theta, \theta_0)} + \frac{C\log b}{I(\theta, \theta_0)}(1+o(1)).
\end{equation}
Note that a similar argument can be found in \cite{wang2015large}.
\end{proof}

Next we will establish a few Lemmas useful for proving theorem \ref{CADD1} for sequential detection procedures. Define a measure $\mathbb{Q}$ on $(\mathcal{X}^{\infty}, \mathcal{B}^{\infty})$ under which the probability density of $X_i$ conditional on $\mathcal{F}_{i-1}$ is $f_{\hat{\theta}_{i-1}}$. Then for any event $A \in \mathcal{F}_{i}$, we have that $\mathbb{Q}(A) = \int_A \Lambda_i d \mathbb{P}_{\infty}$. The following lemma shows that the restriction of $\mathbb{Q}$ to $\mathcal{F}_i$ is well defined. 
\begin{Lemma}
Let $\mathbb{Q}_i$ be the restriction of $\mathbb{Q}$ to $\mathcal{F}_i$. Then for any $A \in \mathcal{F}_k$ and any $i\geq k$, $\mathbb{Q}_i(A) = \mathbb{Q}_k(A)$.
\label{restrictionlemma}
\end{Lemma}

\begin{proof}[Proof of Lemma \ref{typeIerrorbound}]
To bound the term $\mathbb{P}_{\infty}(\tau(b) < \infty)$, we need take advantage of the martingale property of $\Lambda_t$ in (\ref{onestat}). The major technique is the combination of change of measure and Wald's likelihood ratio identity \cite{siegmund1985sequential}. \yc{The proofs are a combination of the results in \cite{lai2004likelihood} and \cite{lorden2005nonanticipating} and the reader can find a complete proof in \cite{lai2004likelihood}. For purpose of completeness we copy those proofs here.}

Define the $L_i = d\mathbb{P}_i/d\mathbb{Q}_i$ as the Radon-Nikodym derivative, where $\mathbb{P}_i$ and $\mathbb{Q}_i$ are the restriction of $\mathbb{P}_{\infty}$ and $\mathbb{Q}$ to $\mathcal{F}_i$, respectively. Then we have that $L_i = (\Lambda_i)^{-1}$ for any $i\geq 1$ (note that $\Lambda_i$ is defined in (\ref{onestat})).
Combining the Lemma \ref{restrictionlemma} and the Wald's likelihood ratio identity, we have that 
\begin{equation}
\mathbb{P}_{\infty}(A \cap \{\tau(b) < \infty\}) = \mathbb{E}_Q \left[\mathbb{I}(\{\tau(b) < \infty\}) \cdot L_{\tau(b)} \right], \forall A\in \mathcal{F}_{\tau(b)},
\label{waldidentity}
\end{equation}
where $\mathbb{I}(E)$ is an indicator function that is equal to 1 for any $\omega \in E$ and is equal to $0$ otherwise.
By the definition of $\tau(b)$ we have that $L_{\tau(b)} \leq \exp(-b)$. Taking $A=\mathcal{X}^{\infty}$ in (\ref{waldidentity})  we prove that $\mathbb{P}_{\infty}(\tau(b)<\infty) \leq \exp(-b)$.
\end{proof}

\begin{proof}[Proof of Corollary \ref{cor:optimalityHP}]
Using (5.180) and (5.188) in \cite{tartakovsky2014sequential}, which are about asymptotic performance of open-ended tests. Since our problem is a special case of the problem in \cite{tartakovsky2014sequential}, we can obtain
\[
\inf_{T\in C(\alpha)}\mathbb{E}_{\theta,0}[T] =  \frac{\log \alpha}{I(\theta, \theta_0)} + \frac{\log (\log(1/\alpha))}{2I(\theta, \theta_0)}(1+o(1)).
\]
Combing the above result and the right-hand side of (\ref{logregretresult}), we prove the corollary.
\end{proof}

\begin{proof}[Proof of Theorem \ref{CADD1}]
From (\ref{ASRandACM}), we have that for any $\nu\geq 1$, 
\[
\mathbb{E}_{\theta, \nu}[T_{ASR}(b)-\nu \mid T_{ASR}(b)>\nu] \leq \mathbb{E}_{\theta, \nu}[T_{ACM}(b)-\nu \mid T_{ACM}(b)>\nu].
\]
Therefore, to prove the theorem using Theorem \ref{maintheorem}, it suffices to show that 
\[
\sup_{\nu\geq 0} \mathbb{E}_{\theta, \nu}[T_{ACM}(b)-\nu \mid T_{ACM}(b)>\nu] \leq \mathbb{E}_{\theta,0}[\tau(b)].
\]
Using an argument similar to the remarks in \cite{lorden2005nonanticipating}, we have that the supreme of detection delay over all change locations is achieved by the case when change occurs at the first instance, \begin{equation}
\sup_{\nu\geq 0} \mathbb{E}_{\theta, \nu}[T_{ACM}(b)-\nu \mid T_{ACM}(b)>\nu] = \mathbb{E}_{\theta,0}[T_{ACM}(b)].
\label{simplification}
\end{equation}
\yc{This is a slight modification (a small change on the subscripts) of the remarks in \cite{lorden2005nonanticipating} but for the purpose of completeness and clearness we write the details in the following}. 
Notice that since $\theta_0$ is known, for any $j \geq 1$, the distribution of $\{ \max_{j+1\leq k\leq t}\Lambda_{k,t} \}_{t=j+1}^{\infty}$ under $\mathbb{P}_{\theta, j}$ conditional on $\mathcal{F}_{j}$ is the same as the distribution of $\{\max_{1\leq k\leq t} \Lambda_{k,t}\}_{t=1}^{\infty}$ under $\mathbb{P}_{\theta, 0}$. Below, we use a renewal property of the ACM procedure. Define 
\[
T_{ACM}^{(j)}(b) = \inf \{t>j: \max_{j+1\leq k\leq t} \log \Lambda_{k,t} >b \}.
\]
Then we have that $\mathbb{E}_{\theta,0}[T_{ACM}(b)]= \mathbb{E}_{\theta, j}[T_{ACM}^{(j)}(b)-j \mid T_{ACM}^{(j)}(b) > j]$. However, $\max_{1\leq k \leq t} \log \Lambda_{k,t} \geq  \max_{j+1\leq k\leq t} \Lambda_{k,t}$ for any $t>j$. Therefore, $T_{ACM}^{(j)}(b) \geq T_{ACM}(b)$ conditioning on $\{T_{ACM}(b)>j \}$. So that for all $j \geq 1$, 
\[
\mathbb{E}_{\theta, 0}[T_{ACM}(b)] = \mathbb{E}_{\theta, j}[T_{ACM}^{(j)}(b) - j \mid T_{ACM}(b)>j] \geq \mathbb{E}_{\theta, j}[T_{ACM}(b)-j\mid T_{ACM}(b)>j].
\]
Thus, to prove (\ref{simplification}), it suffices to show that $\mathbb{E}_{\theta,0}[T_{ACM}(b)] \leq \mathbb{E}_{\theta,0}[\tau(b)]$. To show this, define $\tau(b)^{(t)}$ as the new stopping time that applies the \yc{one-sided }sequential hypothesis testing procedure $\tau(b)$ to data $\{X_i\}_{i=t}^{\infty}$. Then we have that in fact $T_{ACM}(b) = \min_{t\geq 1} \{\tau(b)^{(t)}+t-1 \}$, this relationship was developed in \cite{lorden1971procedures}. Thus, $T_{ACM}(b) \leq \tau(b)^{(1)} + 1 - 1 = \tau(b)$, and $\mathbb{E}_{\theta,0}[T_{ACM}(b)] \leq \mathbb{E}_{\theta,0}[\tau(b)]$. 
\end{proof}

\begin{proof}[Proof of Lemma \ref{ARL1}]
\yc{This is a classic result proved by using the martingale property and the proof routine can be found in many textbooks such as \cite{tartakovsky2014sequential}}. First, rewrite $T_{ASR}(b)$ as
\[
T_{ASR}(b) = \inf\left\{t\geq 1: \log\left(\sum_{k=1}^t \Lambda_{k,t}\right) >b \right\}.
\]
Next, since
\begin{equation}
\log\left(\sum_{k=1}^t \Lambda_{k,t}\right) > \log\left(\max_{1\leq k\leq t} \Lambda_{k,t}\right) = \max_{1\leq k\leq t} \log \Lambda_{k,t},
\label{ASRandACM}
\end{equation}
we have $\mathbb{E}_{\infty}[T_{ACM}(b)]\geq \mathbb{E}_{\infty}[T_{ASR}(b)]$. So it suffices to show that $\mathbb{E}_{\infty}[T_{ASR}(b)]\geq \gamma$, if $b\geq \log \gamma$. 
Define $R_t = \sum_{k=1}^t \Lambda_{k,t}$. Direct computation shows that 
\begin{equation}
\begin{split}
\mathbb{E}_{\infty}[R_t \mid \mathcal{F}_{t-1}] =& \mathbb{E}_{\infty}\left[\Lambda_{t,t} + \sum_{k=1}^{t-1} \Lambda_{k,t} \mid \mathcal{F}_{t-1} \right] \\ 
=& \mathbb{E}_{\infty} \left[1 + \sum_{k=1}^{t-1} \Lambda_{k,t-1} \cdot \log \frac{f_{\hat{\theta}_{t-1}}(X_t)}{f_{\theta_0}(X_t)} \mid \mathcal{F}_{t-1} \right]\\ 
=& 1 + \sum_{k=1}^{t-1} \Lambda_{k,t-1} \cdot \mathbb{E}_{\infty}\left[\log \frac{f_{\hat{\theta}_{t-1}}(X_t)}{f_{\theta_0}(X_t)} \mid \mathcal{F}_{t-1} \right] \\
=& 1 + R_{t-1}.
\end{split} \nonumber
\end{equation}
Therefore, $\{R_t -t \}_{t\geq 1}$ is a $(\mathbb{P}_{\infty}, \mathcal{F}_t)$-martingale with zero mean. Suppose that $\mathbb{E}_{\infty}[T_{ASR}(b)] < \infty$ (otherwise the statement of proposition is trivial), then we have that 
\begin{equation}
\sum_{t=1}^{\infty} \mathbb{P}_{\infty}(T_{ASR}(b) \geq t) < \infty.
\label{temp1}
\end{equation}
(\ref{temp1}) leads to the fact that $\mathbb{P}_{\infty}(T_{ASR}(b)) \geq t = o(t^{-1})$ and the fact that $0\leq R_t \leq \exp(b)$ conditioning on the event $\{T_{ASR}(b) > t\}$,
we have that 
\[
\liminf_{t \rightarrow \infty} \int_{\{T_{ASR}(b)>t\}} |R_t-t| d\mathbb{P}_{\infty} \leq \liminf_{t\rightarrow \infty} ~ (\exp(b)+t)\mathbb{P}_{\infty}(T_{ASR}(b)\geq t) = 0. 
\]
Therefore, we can apply the optional stopping theorem for martingales, to obtain that $\mathbb{E}_{\infty}[R_{T_{ASR}(b)}] = \mathbb{E}_{\infty}[T_{ASR}(b)]$. By the definition of $T_{ASR}(b)$, $R_{T_{ASR}(b)} > \exp(b)$ we have that $\mathbb{E}_{\infty}[T_{ASR}(b)] > \exp(b)$. Therefore, if $b\geq \log \gamma$, we have that $\mathbb{E}_{\infty}[T_{ACM}(b)] \geq \mathbb{E}_{\infty}[T_{ASR}(b)] \geq \gamma$. 
\end{proof}

\begin{proof}[Proof of Corollary \ref{optimality_2}]
Our Theorem 1 and the remarks in \cite{siegmund2008minimax} show that the minimum worst-case detection delay, given a fixed ARL level $\gamma$, is given by
\begin{equation}
\inf_{T(b)\in S(\gamma)}\sup_{\nu\geq 1} \mathbb{E}_{\theta, \nu}[T(b)-\nu+1 \mid T(b)\geq\nu] = \frac{\log \gamma}{I(\theta, \theta_0)} + \frac{d\log \log \gamma}{2I(\theta, \theta_0)}(1+o(1)).
\label{eq:worstEDD}
\end{equation}
It can be shown that the infimum is attained by choosing $T(b)$ as a weighted Shiryayev-Roberts detection procedure, with a careful choice of the weight over the parameter space $\Theta$. 
Combing (\ref{eq:worstEDD}) with the right-hand side of (\ref{logregretresult}), we prove the corollary.
\end{proof}

The following derivation borrows ideas from \citep{azoury2001relative}. First, we derive concise forms of the two terms in the definition of $R_t$ in (\ref{regret}).
\begin{Lemma}
Assume that $X_1, X_2, \ldots$ are i.i.d. random variables with density function $f_{\theta}(x)$, and assume decreasing step-size $\eta_i = 1/i$ in Algorithm \ref{alg1}. Given $\{\hat{\theta}_i\}_{i\geq 1}, \{\hat{\mu}_i\}_{i\geq 1}$ generated by Algorithm \ref{alg1}. If $\hat{\theta}_i = \tilde{\theta}_i$ for any $i \geq 1$, then for any null distribution parameter $\theta_0 \in \Theta$ and $t\geq 1$, 
\begin{equation}
 \sum_{i=1}^t \{-\log f_{\hat{\theta}_{i-1}}(X_i)\} = \sum_{i=1}^t i B_{\Phi^*}(\hat{\mu}_i, \hat{\mu}_{i-1}) - t\Phi^*(\hat{\mu}_t). 
\label{eq:regret_term1}
\end{equation}
Moreover, for any $t\geq 1$, 
\begin{equation}
\inf_{\tilde{\theta} \in \Theta} \sum_{i=1}^t \{-\log f_{\tilde{\theta}}(X_i)\} = -t\Phi^*(\hat{\mu}),
\label{eq:regret_term2}
\end{equation}
where $\hat{\mu} = (1/t)\cdot \sum_{i=1}^t \phi(X_i)$.
\label{lem:regret_term}
\end{Lemma}
By subtracting the expressions in (\ref{eq:regret_term1}) and (\ref{eq:regret_term2}), we obtain the following result which shows that the regret can be represented by a weighted sum of the Bregman divergences between two consecutive estimators.

\begin{proof}[Proof of Lemma \ref{lem:regret_term}]
By the definition of the Legendre-Fenchel dual function we have that
$
\Phi^*(\mu) = \theta^\intercal \mu - \Phi(\theta)
$ 
for any $\theta \in \Theta$. 
By this definition, and choosing $\eta_i = 1/i$, 
we have that for any $i\geq 1$
\begin{equation}
\begin{split}
&-\log f_{\hat{\theta}_{i-1}}(X_i) 
= \Phi(\hat{\theta}_{i-1}) - \hat{\theta}_{i-1}^\intercal \phi(X_i) 
= \hat{\theta}_{i-1}^{\intercal}(\hat{\mu}_{t-1} - \phi(X_i)) - \Phi^*(\hat{\mu}_{i-1}) = \frac{1}{\eta_i}\hat{\theta}_{i-1}^{\intercal} (\hat{\mu}_{i-1} - \hat{\mu}_{i}) - \Phi^*(\hat{\mu}_{i-1})  \\
&=  \frac{1}{\eta_i} (\Phi^*(\hat{\mu}_i) - \Phi^*(\hat{\mu}_{i-1})) - \hat{\theta}_{i-1}^{\intercal} (\hat{\mu}_{i} - \hat{\mu}_{i-1}) - \frac{1}{\eta_i}\Phi^*(\hat{\mu}_i) + \left(\frac{1}{\eta_i}-1\right) \Phi^*(\hat{\mu}_{i-1}) \\
&=  \frac{1}{\eta_i} B_{\Phi^*}(\hat{\mu}_i, \hat{\mu}_{i-1}) + \frac{1}{\eta_{i-1}} \Phi^*(\hat{\mu}_{i-1}) - \frac1{\eta_i} \Phi^*(\hat{\mu}_i), 
\end{split}
\label{eq:lemma3eq1}
\end{equation}
where we use the update rule in Line 6 of Algorithm \ref{alg1} and the assumption $\hat{\theta}_i = \tilde{\theta}_i$ to have the third equation. We define $1/\eta_0 = 0$ in the last equation. Now summing the terms in (\ref{eq:lemma3eq1}), where the second term form a telescopic series, over $i$ from $1$ to $t$, we have that 
\begin{equation}
\begin{split}
\sum_{i=1}^t \{-\log f_{\hat{\theta}_{i-1}}(X_i)\} =& \sum_{i=1}^t \frac{1}{\eta_i} B_{\Phi^*}(\hat{\mu}_i, \hat{\mu}_{i-1}) + \frac{1}{\eta_0} \Phi^*(\hat{\mu}_0) - \frac{1}{\eta_t} \Phi^*(\hat{\mu}_t) \\
=& \sum_{i=1}^t \frac{1}{\eta_i} B_{\Phi^*}(\hat{\mu}_i, \hat{\mu}_{i-1}) - t\Phi^*(\hat{\mu}_t). 
\end{split} \nonumber
\end{equation}
Moreover, from the definition we have that 
\[
\sum_{i=1}^t \{-\log f_{\theta}(X_i)\} = \sum_{i=1}^t \left[ \Phi(\theta) - \theta^\intercal \phi(X_i) \right].
\] 
Taking the first derivative of $\sum_{i=1}^t \{-\log f_{\theta}(X_i)\}$ with respect to $\theta$ and setting it to $0$, we find $\hat{\mu}$, the stationary point, given by
\[
\hat{\mu} = \nabla \Phi(\theta) = \frac{1}{t} \sum_{i=1}^t \phi(X_i). 
\]
Similarly, using the expression of the dual function, and plugging $\hat{\mu}$ back into the equation, we have that 
\begin{equation}\nonumber
\inf_{\tilde{\theta} \in \Theta} \sum_{i=1}^t \{-\log f_{\tilde{\theta}}(X_i)\} =t \cdot \theta^{\intercal} \hat{\mu} - t\Phi^*(\hat{\mu}) - \sum_{i=1}^t \theta^{\intercal} \phi(X_i) = -t\Phi^*(\hat{\mu}).  
\end{equation}
\end{proof}

\begin{proof}[Proof of Theorem \ref{thm:regret}]

By choosing the step-size $\eta_i = 1/i$ for any $i\geq 1$ in Algorithm \ref{alg1}, and assuming $\hat{\theta}_i = \tilde{\theta}_i$ for any $i\geq 1$, we have by induction that 
\[
\hat{\mu}_t = \frac{1}{t} \sum_{i=1}^t \phi(X_i) = \hat{\mu}. 
\]
Subtracting (\ref{eq:regret_term1}) by (\ref{eq:regret_term2}), we obtain
\[
\begin{split}
\mathcal R_t &= \sum_{i=1}^t \{-\log f_{\hat{\theta}_{i-1}}(X_i)\} - \inf_{\tilde{\theta} \in \Theta} \sum_{i=1}^t \{-\log f_{\tilde{\theta}}(X_i)\}\\
& = \sum_{i=1}^t i B_{\Phi^*}(\hat{\mu}_i, \hat{\mu}_{i-1}) - t\Phi^*(\hat{\mu}_t) + t\Phi^*(\hat{\mu})\\
& =  \sum_{i=1}^t i B_{\Phi^*}(\hat{\mu}_i, \hat{\mu}_{i-1})\\
& =  \sum_{i=1}^t i [\Phi^*(\hat{\mu}_i)-\Phi^*(\hat{\mu}_{i-1})-\langle\nabla \Phi^*(\hat{\mu}_{i-1}), \hat{\mu}_i-\hat{\mu}_{i-1} \rangle]\\
& = \frac{1}{2} \sum_{i=1}^t i \cdot (\hat{\mu}_i-\hat{\mu}_{i-1})^\intercal [\nabla^2 \Phi^*(\tilde{\mu}_i)](\hat{\mu}_i-\hat{\mu}_{i-1}).
\end{split}
\]
The final equality is obtained by Taylor expansion. 
\end{proof}

\end{document}